\documentclass[12pt]{article}
\usepackage{amsfonts}
\usepackage{amssymb}
\usepackage{polski}

\input{tcilatex}
\begin{document}

\begin{center}
\textbf{Beyond Erd\H{o}s-Kunen-Mauldin:}

\textbf{Singular sets with shift-compactness properties}

\textbf{by}

\textbf{H. I. Miller\dag \footnote{%
It is with regret that we announce that the first author Harry I. Miller
(1939-2018) died on 16 Dec 2018, whilst finishing this paper. His last
communication ended with: `I only wish I was 10 years younger (make that 15)
so I could contribute (and learn) more.\textquotedblright }, L. Miller-Van
Wieren and A. J. Ostaszewski.}

\bigskip

\bigskip
\end{center}

\noindent \textbf{Abstract.} The Kestelman-Borwein-Ditor Theorem asserts
that a non-negligible subset of $\mathbb{R}$ which is Baire (=has the Baire
property, BP) or measurable is shift-compact: it contains some subsequence
of any null sequence to within translation by an element of the set.
Effective proofs are recognized to yield (i) analogous category and
Haar-measure metrizable generalizations for Baire groups and locally compact
groups respectively, and (ii) permit under $V=L$ construction of co-analytic
shift-compact subsets of $\mathbb{R}$ with singular properties, e.g. being
concentrated on $\mathbb{Q}$, the rationals.

\bigskip


\bigskip

\noindent \textbf{Classification}: 26A03, 04A15, 02K20, 39B62.

\bigskip

\section{Introduction}

This paper is a sequel to [MilO] where two of the present authors studied
shift-compactness (below), a compactness-like embedding property arising
from infinite combinatorics in $\mathbb{R}$, from two points of view:\
topological (group action yielding dual ways of embedding, and so two ways
of asserting the property), and combinatorial (effective embeddings,
employing completeness of $\mathbb{R}$, and limitations exemplified by
`counter-examples'). Here we return to both these themes, motivated
principally by the effectiveness theme. First, we show that effectiveness
allows completeness to be replaced by category: a semitopological Baire
group $X$ will suffice (definitions in \S 2). Secondly, effectiveness
enables `counter-examples' to gain `good topological character': they may be 
\textit{co-analytic} under G\"{o}del's Axiom of Constructibility $V=L$.

In its most useful form and in its simplest context, that of $\mathbb{R}$,
the property of \textit{shift-compactness} of a subset $T$ asserts that some
subsequence $\{z_{m}\}_{m\in \mathbb{M}}$ (for an infinite $\mathbb{M}$ $%
\mathbb{\subseteq N}$) of any given \textit{null sequence} $z_{n}\rightarrow
0$ may be embedded in $T$ under the action of translation. This embedding
idea can be traced back to Banach [Ban, Ch. I, Th. 4], but its explicit
development goes back to Kestelman [Kes1,2] and to Borwein-Ditor [BorD].
Here, the classically familiar non-negligible sets, both the the Baire
non-meagre and the measurable non-null sets, have this property. That is
precisely the content of the Kestelman-Borwein-Ditor Theorem, KBD. Because
of this, it is often possible to unify category and measure arguments, and
so to bring unity to several areas of classical analysis, such as the
automatic continuity results in the theory of functional equations (the
theorems of Ostrowski and Banach-Mehdi concerning the familiar Cauchy
equation, cf. [BinO5], that of Bernstein-Doetsch concerning mid-point convex
functions, cf. [BinO6]), and fundamental results in the theory of regularly
varying (RV) functions (for instance, Karamata's Uniform Convergence Theorem
-- see [BinGT], or Kendall's Theorem, which characterizes RV sequentially,
cf. [BinO11]).

The broader context is that of groups $G$ with some appropriate topological
structure acting on metrizable spaces $X$, the embeddings being provided by
group action (isometries, or more generally homeomorphisms), including that
of a group $X$ acting on itself by translation. So the null sequences now
converge either to the identity map on $X$ or to the neutral element of the
group $1_{X}$. Here the category argument can assume primacy, since it
subsumes the measure analogue, at least in the locally-compact context
provided by Haar measure, by passage to the Haar density topology (under
which null sets become meagre, as first observed by Haupt and Pauc [HauP],
cf. [Kec, 17.47(iii)] and [BinO7, \S\ 2]). Shift-compactness under group
action implies the celebrated `Open Mapping Principle' due to Effros [Eff],
cf. [Ost1,3].

A further key to success in unifying several areas of analysis is the
Steinhaus-Weil Interior Point Theorem [Ste], [Wei], here regarded as
including the Piccard-Pettis Theorem [Pic], [Pet] (since it is true both for
category and measure), that the neutral element is an interior point of $%
AA^{-1}$ for $A$ non-negligible (Baire/measurable). In fact, the theorem
follows from shift-compactness of $A$: see [BinO5], [BinO9], and Theorem
3(ix) below.

To go beyond the Haar context of Polish groups, one needs to abandon measure
invariance, which is prescribed for \textit{all} measurable sets and \textit{%
all} translations. On the measure side an abelian setting is usually (though
not exclusively) preferable and, referring to the family of probability
measures, one needs the \textit{Haar-null} sets of Christensen [Chr1,2],
where \textit{one} particular set remains null under one corresponding
(probability) measure and\textit{\ all} translations; more generally, one
may make do with the (left) Haar-null sets of Solecki [Sol] (albeit aided by
a localized notion of amenability). On the category side there are their
relatives: the \textit{Haar-meagre} sets of Darji [Dar], where the
particular set and all its translations have meagre preimages under some 
\textit{one} continuous map with a compact metric domain; for background see
[Jab]. (These are indeed all meagre.)

But, instead, one may fix \textit{one} reference measure and then use only 
\textit{admissible} translations, rather than all translations, and \textit{%
relative quasi-invariance} of measure (preservation of nullity, relative to
the admissible translations). The canonical example here is a Gaussian
measure in a Hilbert space where the admissible translations form the 
\textit{Cameron-Martin space} [Bog], again a Hilbert space, but under a
refinement of the norm. For literature and generalizations, see [BinO10].

The category-measure duality visible above relies on qualitative aspects of
measure theory, rather than quantitative, and it is \textit{refinement
topologies (density topologies) }which clarify the transition: see [BinO7].

The dichotomy of category -- meagre versus non-meagre sets -- has a
corresponding dichotomy (in an abelian group) between shift-compact and
non-shift-compact sets. The latter have recently been named \textit{%
null-finite} [BanJ], by analogy with Haar null and Haar meagre, and indeed 
\textit{universally measurable} null-finite sets are Haar-null (i.e.
non-Haar-null sets are shift-compact, as has been noted independently in
[BanJ, Th. 4.1] and [BinO9, Th. 3]). Likewise, null-finite sets that are 
\textit{universally Baire} (i.e. pre-images under all continuous maps with a
compact metric domain are Baire) are Haar-meagre [BanJ, Th. 3.1]; for
further background see [BanGJSJ]. The universal Baire property first arose
in mathematical logic: see [FenMW].

In \S\ 2 we re-prove KBD in a Baire-space setting. Effective versions are
shown in \S 3 and used later in \S 5. In \S 4 we take up the study of
singular sets, reviewing some recent results and also adding new ones to the
stock of known examples; here they are often constructed by transfinite
induction. We typify in \S 5 Theorem 4 the detailed treatment needed to
upgrade the topological character by reference to just one of the results
reviewed in \S 4, Theorem MM, by applying G\"{o}del's Axiom of
Constructibility $V=L;$ the other relevant examples of Theorem 3 are
relegated to Theorem 4$^{\prime }$ but with a sketched proof. Our treatment
follows in the footsteps of Erd\H{o}s-Kunen-Mauldin [ErdKM], as in our
title, but we take note of the general `black-box' approach recently
advanced by Vidny\'{a}nszky [Vid] (contemporaneous with our own earlier
development, acknowledged in [MilM], for which it was drafted as supporting
material). We close in\S 6 with complements.

\section{Kestelman-Borwein-Ditor Theorem: topological setting}

There are a number of versions of the KBD and so of its proof, which go back
to [Kes1,2], [BorD] -- for an account see [BinO5], [BinO4], [BinO3] and
[MilO]. This section is dedicated to a proof applicable to the context of a
Baire \textit{semitopological group} (defined below), Theorem 2, based on
the proof strategy used in [MilO] to prove KBD in $\mathbb{R}.$ Although on
first inspection it may seem that that proof, in constructing inductively a
sequence of approximations to a translator, uses completeness of $\mathbb{R}$%
, and so is adaptable only to a completely metrizable space, in fact matters
are otherwise. The inductive step is sufficiently typical, i.e. unspecific
to the preceeding step, that it may be applied anywhere in space; so the
Baire theorem will carry through the induction `to the limit' at least
somewhere (and so almost everywhere, according to the Generic Dichotomy
Principle [BinO2]).

We close the section with the statement of another version of the KBD
applicable to \textit{topological groups}, one that is strong enough to
imply the celebrated result of Effros [Eff, \S\ 2] known as the Open Mapping
Principle [Anc], cf. [Ost1,3]. As one would expect this does indeed imply
Theorem 2 when specialized to topological groups (see the Theorem from
[Ost3] at the end of the section)

We first prove in Theorem 1 a special case of our KBD here, and then deduce
the main result as Theorem 2. That is followed by its \textit{Haar-measure}
version, Theorem 2H. The latter closely follows the argument in [MilO, Th.
1M]. However, the present Haar context calls for some extra details.

Below $z_{0}$ is the identity element of the group. Also we recall that a
group is said to be \textit{semitopological} [ArhT] if translation is
continuous -- so an autohomeomorphism. A \textit{space} is \textit{Baire} if
it obeys Baire's Category Theorem; however, a \textit{set} $A$ is \textit{%
Baire} if it has the Baire property, BP. Then we denote by $A^{q}$ the 
\textit{quasi-interior} of $A,$ i.e. the largest open set equal to $A$
modulo a meagre set.

\bigskip

\noindent \textbf{Theorem 1.} \textit{In a Baire semitopological group }$X$%
\textit{, if }$A$\textit{\ is co-meagre and }$\{z_{n}\}_{n\in \mathbb{N}%
}\rightarrow 1_{X}$\textit{\ is a null sequence, then for a dense }$\mathcal{%
G}_{\delta }$\textit{-set of points }$a$\textit{\ in }$A:$ 
\[
\{az_{n}:n=0,1,2,...\}\subseteq A. 
\]

\bigskip

For subsets $A,B$\ of $X$ below we write $AB:=\{ab:a\in A,b\in B\}.$ We
abbreviate neighbourhood to \textit{nhd} and nowhere dense to \textit{nwd}.
We begin with

\bigskip

\noindent \textbf{Lemma 1 (Extension of a separation).} \textit{In a
semitopological group }$X$\textit{, for }$f\in X$\textit{, finite }$%
F\subseteq X,$\textit{\ }$L$\textit{\ nwd, and }$V$\textit{\ non-empty open,
if }%
\[
(FV)\cap L=\emptyset , 
\]%
\textit{then there is a non-empty open }$V^{\prime }\subseteq V$\textit{\
with}%
\[
((F\cup \{f\})V^{\prime })\cap L=\emptyset . 
\]

\noindent \textbf{Proof.} Given $L,f,F$ and open $V,$ as $fU\ $is non-empty
and open, choose a non-empty open $U\subseteq fV$ with $\emptyset =U\cap L.$
Put $V^{\prime }:=f^{-1}U\subseteq V;$ then%
\[
(fV^{\prime })\cap L=U\cap L=\emptyset , 
\]%
and%
\[
(FV^{\prime })\cap L\subseteq (FV)\cap L=\emptyset . 
\]%
Since $(F\cup \{f\})V^{\prime }=(FV^{\prime })\cup (fV^{\prime }),$%
\[
((F\cup \{f\})V^{\prime })\cap L=\emptyset .\qquad \square 
\]

\bigskip

\noindent \textbf{Proof of Theorem 1.} W.l.o.g. the sequence is injective,
so in particular (with $z_{0}=1_{X})$ $Z:=\{z_{i}:i=0,1,2,...\}$ is
infinite. We put $Z_{n}:=\{z_{0},z_{1},...,z_{n}\}.$ For any finite set of
points $F=\{f_{0},f_{1},...,f_{k}\},$ with $f_{0}=z_{0}=1_{X},$ $L$ nwd, and 
$G\ $open and dense in $X,$ for $0\leq i\leq n$ put $F_{i}:=%
\{f_{0},...,f_{i}\}$ (so that $F=F_{n}),$ and let%
\[
W_{L}^{F}(G):=\{x\in G:(\exists \text{ open }W_{x})[x\in W_{x}\text{ and }%
(FW_{x})\cap L=\emptyset ]\}\subseteq G\backslash L. 
\]%
Then $W:=W_{L}^{F}(G)$ is open, since $W_{x}\subseteq W$ for each $x\in W.$
Notice that $V_{x}:=(x^{-1}W_{x})$ is a nhd of $1_{X}$ (since translation is
a homeomorphism) with 
\[
\lbrack (Fx)V_{x}]\cap L=0, 
\]%
i.e. $V_{x}$ generates a nhd of the shifted set $Fx$ disjoint from $L.$

\textit{Claim 1. }The open set $W$ is dense in $G$.

For non-empty open $U\subseteq G,$ define inductively non-empty open sets $%
V_{i}$ with $U\supseteq V_{0}\supseteq ...V_{i-1}\supseteq V_{i}\supseteq
...\supseteq V_{n}$ such that for $0\leq i\leq n$ 
\[
(F_{i}V_{i})\cap L=\emptyset , 
\]%
i.e. for any $x\in V_{i}$ , $V_{x}^{i}=x^{-1}V_{i}$ is a nhd of $1_{X}$ with 
\[
\lbrack (F_{i}x)V_{x}^{i}]\cap L=\emptyset , 
\]%
so providing a uniform nhd for the shifted set $F_{i}x$ disjoint from $L.$

\textit{Basis step}. As $L$ is nwd, choose non-empty open $V_{0}\subseteq U$
with $\emptyset =V_{0}\cap L=(F_{0}V_{0})\cap L,$ as $F_{0}=\{z_{0}\}=%
\{1_{X}\}.$

\textit{Inductive step. }Given $V_{i-1},$ apply Lemma 1 to $L,f_{i},F_{i-1}$
and $V_{i-1}$ to choose a non-empty open $V^{\prime }$ as in the Lemma. Take 
$V_{i}:=V^{\prime };$ then\textit{\ }$V_{i}\subseteq V_{i-1}$ and%
\[
(F_{i}V_{i})\cap L=((F_{i-1}\cup \{f_{i}\})V_{i})\cap L=\emptyset . 
\]%
At the conclusion of the induction, for $x\in V_{n}$ the set $W_{x}:=V_{n}$
gives%
\[
(FW_{x})\cap L=\emptyset , 
\]%
and so $x\in W\cap U,$ proving density of $W$ in $G$. $\square _{\text{%
(claim 1)}}$

\textit{Claim 2.} For $F^{\prime }=F\cup \{f\},$ and $G=W$ above, $%
W_{L}^{F^{\prime }}(G)$ is open and dense in $G.$

This is almost a repeat of the inductive step in the proof of Claim 1. For
non-empty open $U\subseteq G,$ w.l.o.g. we may assume by Claim 1 that $%
U\subseteq W=W_{L}^{F}(G),$ as $W$ is dense open. Consider any $x\in U.$ As $%
x\in W_{L}^{F}(G),$ there is an open nhd $W_{x}$ of $x$ with $(FW_{x})\cap
L=\emptyset .$ W.l.o.g. $W_{x}\subseteq U.$ (Indeed, as $x\in U\cap
W_{x}\subseteq W_{x}$, we have $[F(U\cap W_{x})]\cap L=\emptyset .$)

Apply Lemma 1 to $L,f,F$ and $W_{x}$ to choose a non-empty open $V^{\prime
}\subseteq W_{x}$ with%
\[
(F^{\prime }V^{\prime })\cap L=\emptyset . 
\]%
So $y\in W_{L}^{F^{\prime }}(G),$ for any $y\in V^{\prime }\subseteq
W_{x}\subseteq U$ (with $V^{\prime }$ doing duty for $W_{y}$). $\square _{%
\text{(claim 2)}}$

\bigskip

Now write $X\backslash A=\bigcup\nolimits_{n=0}^{\infty }N_{n}$ with the $%
N_{n}$ increasing and nwd. Put $N_{-1}=\emptyset ,$ $W_{-1}=X,$ and define
inductively dense open sets%
\[
W_{2n}=W_{N_{n-1}}^{Z_{n}}(W_{2n-1}),\qquad
W_{2n+1}=W_{N_{n}}^{Z_{n}}(W_{2n}), 
\]%
so that%
\begin{eqnarray*}
W_{0} &=&W_{\emptyset }^{\{z_{0}\}}(X)=X,\qquad
W_{1}=W_{N_{0}}^{Z_{0}}(X),\qquad W_{2}=W_{N_{0}}^{Z_{1}}(W_{1}), \\
W_{3} &=&W_{N_{1}}^{Z_{1}}(W_{1}),\qquad \text{etc.}
\end{eqnarray*}

Now put%
\[
H=\bigcap\nolimits_{n=0}^{\infty }W_{n}, 
\]%
which, since $X$ is Baire, is dense. Fix $x\in H,$ and $n$; then for $m\geq
n,$ since $x\in W_{2m+1}$ 
\[
xZ_{n}\subseteq xZ_{m}\subseteq X\backslash N_{m} 
\]%
and so%
\[
xZ_{n}\subseteq \bigcap\nolimits_{m=n}^{\infty }X\backslash
N_{m}=X\backslash \bigcup\nolimits_{m=n}^{\infty }N_{m}=A, 
\]%
since the sequence $N_{n}$ is increasing. But $n$ was arbitrary, so%
\[
xZ\subseteq A.\qquad \square 
\]

We deduce as an easy corollary the category version of the KBD:

\bigskip

\noindent \textbf{Theorem 2 (Theorem KBD).} \textit{In a Baire (metric)
semitopological group }$X,$\textit{\ if }$\{z_{n}\}_{n\in \mathbb{N}}$%
\textit{\ is a null sequence, and }$A$\textit{\ is a non-meagre Baire set,
then for some }$a\in A$\textit{\ and some }$n=n(a)$%
\[
\{az_{m}:m>n\}\subseteq A. 
\]%
\textit{In fact, the embedding holds for quasi almost all }$a\in A.$

\bigskip

\noindent \textbf{Proof.} W.l.o.g. $A=A^{q}\backslash M$ with $M\subseteq
A^{q}$ meagre and $A^{q}$ non-empty. As $X\backslash M$ is co-meagre, by Th.
1 there is a\textit{\ }$\mathcal{G}_{\delta }$-set $H$ of points $x$ in $%
X\backslash M$ with 
\[
\{xz_{n}:n=0,1,2,...\}\subseteq X\backslash M. 
\]%
As $A^{q}$ is open, we may choose $a\in H\cap A^{q}\subseteq A^{q}\backslash
M=A.$ As $a=\lim_{n}(az_{n}),$ by continuity of left translation, there is $%
n(a)$ with $\{az_{n}:n>n(a)\}\subseteq A^{q}.$ But then%
\[
\{az_{n}:n>n(a)\}\subseteq A^{q}\cap (X\backslash M)=A. 
\]

The final conclusion holds by the Generic Dichotomy Theorem [BinO2]. $%
\square $

\bigskip

Note the following immediate corollary, concerning measurable groups [Hal, 
\S 62] equipped with a probability measure $\mu $ and a `differentiation
basis' [Bru], giving rise to a density topology $\mathcal{D}_{\mu }$ in the
sense of Martin [Mar]. (The Haar density topology case, using a
differentiation basis provided by [Mue], is discussed in [BinO3, \S 7]; for
background on density topologies see [BinO7].)

\bigskip

\noindent \textbf{Corollary.} \textit{For a topological group }$X$\textit{\
supporting a density topology }$\mathcal{D}_{\mu }$\textit{\ generated by a
measure }$\mu $ \textit{(e.g. Haar measure on a locally compact group): if }$%
A$\textit{\ is co-null and} $\{z_{n}\}_{n\in \mathbb{N}}$\textit{\ is a null
sequence, then for a dense }$\mathcal{G}_{\delta }(\mathcal{D}_{\mu }%
\mathcal{)}$\textit{-set of points }$a$\textit{\ in }$A:$ 
\[
\{az_{n}:n=0,1,2,...\}\subseteq A. 
\]

\noindent \textbf{Proof. }Under the density topology the group is both Baire
[BinO7, Prop 4] and semitopological. As the nwd sets are precisely the $\mu $%
-null sets [BinO7, Th. 7.2], Th. 1 applies. $\square $

\bigskip

The next result, which emerges as more demanding, goes beyond a co-null
setting.

\bigskip

\noindent \textbf{Theorem 2H.}\textit{\ Let }$G$ \textit{\ be a locally
compact metrizable topological group.}

\noindent (i) \textit{For any convergent sequence }$\{x_{n}\}_{n\in \mathbb{N%
}}$\textit{\ with limit }$x_{0}$\textit{\ and any (right) non-null
Haar-measurable set }$T,$\textit{\ there are a left shift }$\theta (x)=cx$%
\textit{\ and an infinite set} $\mathbb{M\subseteq N}$ \textit{such that }$%
\theta (x_{0})\in T$\textit{\ and} 
\[
\theta (x_{m})=cx_{m}\in T\text{ for }m\in \mathbb{M}. 
\]%
\noindent (ii) \textit{Moreover, for }$S$\textit{\ and }$T$\textit{\
density-open with }$Sx_{0}\subseteq T$\textit{\ the shift may be chosen with 
}$c\in S.$

\bigskip

\noindent \textbf{Proof.} Below $|.|$ denotes a right-invariant Haar measure
on $G.$ By the Birkhoff-Kakutani metrization theorem ([Bir], [Kak1], [DieS],
or [Ost2]), we may equip $G$ with a group norm $||x||:=d(1_{G},x)$ for $d$ a
right-invariant metric.

\noindent (i) Let $T\ $be Haar-measurable non-null. By inner regularity of
the measure, we may assume that $T$ is compact and non-null. Applying a left
shift to the sequence $x_{n}$ if necessary, $x_{0}$ is w.l.o.g. a density
point of $T.$ Put $m(0):=0,$ $\theta _{0}:=\mathrm{id.}$

Suppose inductively that $\theta _{n}(x):=c_{n}...c_{1}x$ for some $c_{i}$
with $||c_{i}||\leq 2^{-i}$ for $i\leq n,$ and that an increasing sequence
of integers $m(j)$ for $j\leq n$ has been selected with each $u_{j}:=\theta
_{n}(x_{m(j)})$ a density point of $T$.

For each $\varepsilon =2^{-n}$ we may choose $U$ a finite union of (left)
translates of an open nhd of $1_{G}$ to cover $T$ with the complement $%
E=U\backslash T$ having $|E|<\varepsilon $. Choose open nhds $I_{j}$ with $%
u_{j}\in I_{j}\subseteq U,$ for $0\leq j\leq n.$ Let $\eta :=\min_{0\leq
j\leq n}\{d(u_{j},X\backslash I_{j}),\varepsilon \}.$ Since each $u_{j}$ is
a density point, choose a symmetric open nhd $V$ round $1_{G}$ such that
each $V_{j}:=Vu_{j}\subseteq I_{j}$ has $|V_{j}|=|V|$, $u_{j}\in V_{j}$ and $%
|V_{j}\cap T|\geq (1-\eta )|V|$ for all $j\leq n$ and $|V|<\varepsilon .$
Choose $m=m(n+1)$ with $m>m(n)$ such that $d(x_{m},x_{0})<\eta $ and $%
u_{n+1}:=\theta _{n}(x_{m(n+1)})\in V_{0};$ both are possible as $%
x_{0}=\lim_{m}x_{m}$ and $\theta _{n}(x_{0})=u_{0}\in V_{0}$ and $\theta
_{n} $ is continuous. Choose an open interval $V_{n+1}\subseteq I_{0}$
centered on $u_{n+1}.$

For $j\leq n$ one has $|V_{j}\cap E|<\eta |V|$ as $V_{j}\backslash
E\subseteq U\backslash E\subseteq T,$ and so $|V_{j}\backslash E|>(1-\eta
)|V_{j}|.$ Invoking the Haar Density Theorem ([Mue], [Mar]), let $F$ be a
measure-zero set such that $(V_{j}\backslash E)\backslash F$ is a
density-open subset of $T$ (all its points are density points) for each $j<n$%
.

For any $c,$ note that $cu_{n+1}$ is a density point of $T\cap V_{0}$ iff $c$
is a density point of $T^{\prime }:=(T\cap V_{0})u_{n+1}^{-1},$ as the
measure is right-invariant. Again by the Haar Density Theorem, off a null
subset $N$ of $T^{\prime }$ all its members are density points. In what
follows we ensure that $c\notin N.$

Choose $c_{n+1}\in V\backslash (N\cup (E\cup F)u_{j}^{-1})$ with $%
||c_{n+1}||<\varepsilon $ such that $c_{n+1}u_{j}\in V_{j}\backslash
E_{0}\subseteq T$ and $c_{n+1}u_{j}$ is a density point of $T,$ for each $%
j\leq n+1.$

Set $\theta _{n+1}(x):=c_{n+1}\theta _{n}(x);$ then, for each $j\leq n+1,$ $%
\theta _{n+1}(x_{m(j)})$ is a density point of $T$ in $T.$

Moreover, $s_{n}:=(c_{n}...c_{1})$ converges, to $s$ say, as $%
d(c_{n+1}c_{n}...c_{1},c_{n}...c_{1})=d(c_{n+1},1_{G})$, by right-invariance
of $d,$ and so $\{s_{n}\}_{n\in \mathbb{N}}$ is a Cauchy sequence in a
locally compact nhd of $1_{G}$. Take $\theta (x):=sx;$ then, for each $j,$
as $T$ is compact, $\theta (x_{m(j)})=\lim_{n}\theta _{n}(x_{m(j)})\in T$.
Also $\lim_{j}\theta (x_{m(j)})=\theta (x_{0})\in T,$ as $%
x_{0}=\lim_{m}x_{m}.$

\noindent (ii) This now follows quite easily. Specialize the sequence
arising in the proof above to a null sequence $z_{n}\rightarrow z_{0}=1_{G}$
and replace $T$ by $S$ to obtain $\theta (z_{0})=s1_{G}\in S$ and $sz_{m}\in
S,$ for an infinite set of $m,$ in $\mathbb{M}_{s}$ say.

Returning to a general sequence $x_{n}$ with limit $x_{0},$ put $%
z_{n}:=x_{n}x_{0}^{-1}.$ Then, as before, for some $s\in S$ and some
infinite set $\mathbb{M}_{s},$ one has $sz_{m}\in S$ for $m\in \mathbb{M}%
_{s} $. But then $sz_{m}x_{0}=sx_{m}\in Sx_{0}\subseteq T$ for $m\in \mathbb{%
M}_{s},$ as asserted. $\square $

\bigskip

\noindent \textbf{Remark}. A locally compact metrizable group, being
topologically complete, is completely metrizable [Enge, 4.3.26]. More
generally, it is a theorem of Loy and of Christensen that a topological
Baire group which is analytic is in fact Polish -- see e.g. [TopH, Th.
2.3.6].

\bigskip

We close the section by recalling the promised `strong' version of KBD that
is applicable to topological groups.

\bigskip

\noindent \textbf{Definition }[Ost3],\textbf{\ }cf. [Pet]. For $G$ a
metrizable group, say that the group action $\varphi :G\times X\rightarrow X$
is a \textit{Nikodym }action (or that it has the \textit{Nikodym property})
if for every non-empty open neighbourhood $U\ $of $1_{G}$ and every $x\in X$
the set $Ux=\varphi _{x}(U):=\varphi (x,U)$ contains a non-meagre \textit{%
Baire set}.

\bigskip

\noindent \textbf{Example.} For $G$ a semitopological group acting on
itself: $\varphi _{x}(u)=\varphi (x,u)=ux;$ so $\varphi $ is separately
continuous and $\varphi _{x}$ is an autohomeomorphism. So $Ux$ is open, for
any open $U.$ In particular, if $G$ in Theorem2$\ $is a \textit{topological}
Baire group acting on itself, that action has the Nikodym property, so the
following result implies the conclusion of Theorem 2.

\bigskip

\noindent \textbf{Shift-compactness Theorem} [Ost3]\textbf{.} \textit{For }$%
T $\textit{\ a Baire non-meagre subset of a metric space }$X$\textit{\ and }$%
G$\textit{\ a group, Baire under a right-invariant metric, and with
separately continuous and transitive Nikodym action on }$X$\textit{:}

\textit{for every convergent sequence }$\{x_{n}\}_{n\in \mathbb{N}}$\textit{%
\ with limit }$x_{0}$\textit{\ and any Baire non-meagre }$A\subseteq G$ 
\textit{with }$1_{G}\in A^{q}$ \textit{and} $A^{q}x\cap T^{q}\neq \emptyset
, $\textit{\ there are }$\alpha \in A$\textit{\ and an integer }$N$\textit{\
such that }$\alpha x_{0}\in T$\textit{\ and}%
\[
\{\alpha (x_{n}):n>N\}\subseteq T. 
\]

\section{KBD: effective version}

In this section we give in Theorem 1E an \textit{effective version} of KBD.
Our treatment below of coding overlaps with that of the corresponding
sections of the contemporaneous paper [BinO11].

In what follows, we will need to distinguish between (general) sets of
reals, and `nice' sets which can be defined by a suitable (effective) coding
so that an individual set is coded by a single \textit{real}. For background
here, see e.g. the monograph Kechris [Kec, Ch.V] on the analytical hierarchy
(note [Kec, V.40B] on classical v. effective descriptive set theory), [Rog2,
Part 4] and our recent survey [BinO8]. For a deeper analysis of coding see
[Solo, II.1.1, 25-33]; a minimal amount is in [FenN, \S\ 2, p. 93].

We begin with a short introduction to this topic in the next sub-section (on
preliminaries), which the expert reader can omit. The non-expert reader may
also take `coding on trust', observing the basic case of an open set $%
W\subseteq \mathbb{R}$ which may be coded by first enumerating (effectively)
the rational intervals as $\{I_{n}\}_{n\in \mathbb{N}}$ and then coding $W$
up as the \textit{binary real} which is the indicator function $\mathbf{1}_{%
\mathbb{M}}$ of the subset $\mathbb{M}:=\{n:I_{n}\subseteq W\},$ and thus
omit \S 3.1 We defer further discussion of some of the finer points to the
Appendix of this paper.

\subsection{Preliminaries on coding}

We work in the space $I$ of irrationals, interpreted as the non-recurring 
\textit{binary} sequences $x:\mathbb{N}\rightarrow \{0,1\};$ here $x$ may
also be viewed as the indicator function of a subset of $\mathbb{N}$ and
thereby as a real number \textit{code} for that subset. We may identify $%
x\in I$ with the sequence $\{x_{n}\}_{n\in \mathbb{N}}$, where $x_{n}$
denotes the $n^{\text{th}}$ \textit{projection} defined by $%
x_{n}(m)=x(2^{n}(2m+1)).$ With $x$ viewed as (a code for) a subset of $%
\mathbb{N}$, $x_{n}$ is a code for $x\cap \{1\cdot 2^{n},3\cdot 2^{n},5\cdot
2^{n},$ $...\}.$

The proof of Theorem 3 in \S\ 4 below relies on the ability to refer to
various subsets of the real line in terms of real numbers; in particular, an
open set $G$, closed set $F$, and $\mathcal{G}_{\delta }$-set $H$ may be
coded by $a\in I$ via one of%
\[
G(a):=\bigcup\nolimits_{n\in a}I_{n},\qquad F(a)=I\backslash G(a),\qquad
H(a):=\bigcap\nolimits_{n\in \mathbb{N}}G(a_{n}), 
\]%
where as above $\{I_{n}\}_{n\in \mathbb{N}}$ enumerates (constructively) all
the rational-ended intervals and $a_{n}$ is the $n^{\text{th}}$ projection
of $a$. (Evidently, one must separately code which of the three displayed
equations is to be choosen.) Coding clarifies in what form the property of
`membership in $G(a)$', or $F(a)$ etc., is expressible as an (arithmetic)
predicate in the language of set theory (below); indeed, 
\[
x\in G(a)\text{ iff }(\exists n\in \mathbb{N})[(n\in a)\&(x\in I_{n})],\quad
x\in F(a)\text{ iff }(\forall n\in \mathbb{N})[(n\in a)\&(x\notin I_{n})]. 
\]%
Both predicates here uses an \textit{arithmetic} quantifier (ranging over $%
\mathbb{N}$, the type 0 objects ) while its matrix (the part without
quantifiers, delimited here by square brackets) refers to elementary
relations. The first is said to be $\Sigma _{1}^{0}(a):$ this identifies a
single \textit{existential} quantification over type 0 objects, and the
presence of $a$; its complement is $\Pi _{1}^{0}(a),$ with $\Pi $ for the 
\textit{universal} quantifier. One may supress the explicit mention of $a$
by use of \textit{bold-face} symbols $\mathbf{\Sigma }_{1}^{0}$ and $\mathbf{%
\Pi }_{1}^{0},$ which imply the need for a parameter. By contrast,%
\[
x\in H(a)\text{ iff }(\forall n\in \mathbb{N})(\exists m\in \mathbb{N}%
)[(m\in a_{n})\&(x\in I_{m})], 
\]%
which is $\mathbf{\Pi }_{2}^{0}$ because there is a universal quantifier
leading the alternating \textit{pair} of quantifiers. Similar conventions
govern \textit{analytic} quantifiers (ranging over $\mathbb{R}$, the type 1
objects): the superscript here changes to a 1.

Codes in $\mathbb{R}$ known as `notations' are also needed for the countable
ordinals. This is somewhat tedious, so omitted here. (One may start with the
indicator function of $\mathbb{N}$ as a code for $\omega $ as an order type.)

We make use of the \textit{language of set theory}, $LST:$ its (first-order)
formulae, needed here and again in \S 4, are written using: free variables,
symbols denoting constants, the relation of membership, the usual
connectives, negation, and quantifiers ranging over sets. This enables us to
recall the constructible hierarchy $\langle L_{\alpha }:\alpha \in On\rangle 
$, in which, for ordinal $\alpha ,$ the sets $L_{\alpha }$ are obtained by
iterating transfinitely the operation which defines $L_{\beta +1}$ as the
family of those subsets of $L_{\beta }$ that are definable by the
first-order formulas of $LST.\ $Here they are allowed to refer to a finite
string of elements of $L_{\beta }$ and all their quantifiers range over $%
L_{\beta }$ -- see e.g. [Sac, 9.2.III], [Dev], or [BinO8, \S 2]. The class $%
L:=\tbigcup \{L_{\alpha }:\alpha \in On\}$ comprising all the constructible
sets has a \textit{canonical well-ordering} $<_{L}$(defined by transfinite
induction using an effective listing of all predicates).

\subsection{KBD: a version effective in the codes}

We develop a version of KBD\ suitable for work in $\mathbb{R}.$ We begin by
demonstrating that, for a null sequence $\{z_{n}\}_{n\in \mathbb{N}}$ coded
by a $z\in I$ (with $z_{n}$ as its $n^{\text{th}}$ projection) and for a
target $\mathcal{G}_{\delta }$-set coded by $s\in I,$ the relevant
translator $t$ may be constructed effectively (recursively) in $z$ and $s$.
This guarantees that when the $\mathcal{G}_{\delta }$-set has code $s$ in $%
L_{\alpha }$, then such a translator is in $L_{\alpha +\omega }$ (for $%
L_{\alpha }$ point-definable, as in the preamble in \S 5 to the proof of
Theorem 4). The corresponding $\mathcal{G}_{\delta }$-sets/codes form the
family $\mathcal{G}_{\alpha }$ defined below. (Later on we will also require
the sets $L_{\alpha }$ to be models of the axiom system $\mathrm{ZF}^{-}$($%
\mathrm{ZF}$ less the $\mathrm{Power}$ \textrm{Axiom}); we note that if $%
L_{\alpha }$ is point-definable, then so is $L_{\alpha +\omega }$ -- this is
proved in [EngMS].)

\bigskip

\noindent \textbf{Definitions. }1. Following [MilO], say that the group of
translations $Tr(\mathbb{R}^{d})$ \textit{strongly }$L_{\alpha }$\textit{%
-separates points from a family }$\mathcal{F}$ \textit{of closed nowhere
dense sets} in $\mathbb{R}^{d}$ if for each $p\in L_{\alpha }$ and $F\in 
\mathcal{F}$ and arbitrarily small $q\in \mathbb{Q}_{+}$ there is $%
H\subseteq (-q,q)$ with code in $L_{\alpha }$ such that $h_{c}(p):=p+c\notin
F$ for every $c\in H.$

2. Denote by $\mathcal{G}_{\alpha }$ the family of sets $G$ open (in $%
\mathbb{R}^{d}$) with\textit{\ }$\mathbb{Q}\subseteq G$ possessing a code in 
$L_{\alpha },$ and by $\mathcal{F}_{\alpha }$ the complements of sets in $%
\mathcal{G}_{\alpha }$\textit{.}

3. $B_{\varepsilon }(x)$ denotes the ball centred at $x$ of radius $%
\varepsilon .$

\bigskip

\noindent \textbf{Proposition 1\ (Strong Separation}, cf. [MilO, Prop. 1]%
\textbf{).} \textit{For }$\mathbb{R}^{d}$\textit{\ and }$Tr(\mathbb{R}^{d})$%
\ \textit{both equipped with the Euclidean topology, the group }$Tr(\mathbb{R%
}^{d})$ \textit{strongly }$L_{\alpha }$\textit{-separates points of }$%
L_{\alpha }$\textit{\ and the closed nowhere dense sets of }$\mathcal{F}%
_{\alpha }$\textit{.}

\bigskip

\noindent \textbf{Proof. }Let $q_{i}$ be an effective enumeration of $%
\mathbb{Q}$. For $0<q\in \mathbb{Q}$, if $p\in L_{\alpha }$ and $F=\mathbb{R}%
\backslash G$ with $G$ open and coded in $L_{\alpha },$ choose the first $%
q_{i}\in B_{q}(p)$ and thereafter the first pair $\langle
q_{L(i)},q_{R(i)}\rangle $ with $q_{L(i)}<q_{i}<q_{R(i)}$ such that $%
I:=(q_{L(i)},q_{R(i)})\subseteq G\cap B_{q}(p).$ Then $H:=I-p\subseteq
(-q,q) $ has code in $L_{\alpha },$ and, for $c=m-p\in H$ with $m\in I,$ $%
|c|=|m-p|<q$ and $p+c=m\in G=\mathbb{R}\backslash F$. $\square $

\bigskip

The next result follows, as it is an inductive construction applying Prop. 1
at each inductive stage.

\bigskip

\noindent \textbf{Proposition 2\ (Finitary Euclidean Strong Separation}, cf.
[MilO, Prop. 2])\textbf{.} \textit{With }$F\in \mathcal{F}_{\alpha }$\textit{%
\ as above, let }$U$\textit{\ be Euclidean open with code in }$L_{\alpha }$ 
\textit{and }$u_{i}\in U\ $\textit{for }$i\leq n$ \textit{with }$u_{i}\in
L_{\alpha }.$ \textit{Then, for each }$\varepsilon >0$\textit{, in }$%
B_{\varepsilon }(0)$\textit{\ there is a neighbourhood of }$c$\textit{%
-shifts }$x\rightarrow x+c$\textit{\ with code in }$L_{\alpha }$\textit{\
such that }$u_{i}+c\in U$\textit{\ and }$u_{i}+c\notin F$\textit{\ for each }%
$i\leq n.$

\bigskip

\noindent \textbf{Proof. }Proceed exactly as in [MilO, Prop. 2], using Prop.
1 here in place of Prop. 1 there. $\square $

\bigskip

\noindent \textbf{Theorem 1E} (cf. [MilO, Th. 1E])\textbf{.}\textit{\ For
the real line under the Euclidean topology, given }$y\in \mathbb{N}^{\mathbb{%
N}}\cap L_{\alpha }$ \textit{coding a convergent sequence }$\{y_{n}\}_{n\in 
\mathbb{N}}\rightarrow y_{0},$ \textit{and }$x\in \mathbb{N}^{\mathbb{N}%
}\cap L_{\alpha }$ \textit{such that the set }$T=\bigcap G_{n}\ $\textit{%
with }$G_{n}$ \textit{coded by }$x_{n}$ \textit{and }$G_{n}\in \mathcal{G}%
_{\alpha },$ \textit{there are a }$c$\textit{-shift }$h(x)=x+c$\textit{\ and
an integer }$M$ \textit{such that }$h(x_{0})\in T$\textit{\ and} 
\[
h(y_{m})=y_{m}+c\in T\text{ for }m>M, 
\]%
\textit{and }$c$\textit{\ has a code in }$L_{\alpha +\omega }.$

\bigskip

\noindent \textbf{Proof. }Again proceed exactly as in [MilO, Th. 1E], using
Prop. 2 here in place of Prop. 2 there. $\square $

\section{Singular sets}

In [MilO] the first and third authors studied extensions of the
Kestelman-Borwein-Ditor theorem from the perspective of group action, on the
one hand, and certain limitations (exemplified by `singular' sets) of the
infinite combinatorics involved, on the other. The latter included `wild'
2-place-function actions in place of group actions [MilO, Th. 8], and
examples of non-shift-compactness (existence, for a given closed nwd set $A,$
of a monotonic null sequence with $\{d_{n}\}_{n\in \mathbb{N}}$ so that for
each $x,$ $x+d_{n}\notin A$ infinitely often, and an example, under the
Axiom of Choice, of a non-measurable $A$ with $x+(1/n)\notin A$ for all $n).$

Here, in similar spirit, we offer further examples of singular behaviour. We
recall a particular result needed quite soon. (Below $d(A):=A-A=\{a-a^{%
\prime }:a,a^{\prime }\in A\};$ for $\mathcal{E}mb$ see the Definitions
hereunder.) In the theorem below, the first assumption needed for (i) may be
regarded as a topological variant of Martin's Axiom $(MA):$ see [MilO],
[BinO8, \S 6b]; in (ii) $W$ is concentrated on $\mathbb{Q}$ [Rog1, \S 2.3],
and such sets are of (strong) measure zero: see \S 6.3.

\bigskip

\noindent \textbf{Theorem MO }([MilO, Th. 9]). \textit{Assume that the union
of fewer than }$\mathfrak{c}$\textit{\ many sets that are meagre and null is
itself meagre and null, then there exists }$A\in \mathcal{E}mb$ \textit{with}
$d(A)=\mathbb{R}$ \textit{such that:}

\noindent (i)\textit{\ }$(\forall x\in A)$\textit{\ }$x+(1/n)\notin A$%
\textit{\ for all }$n$\textit{\ with at most one exception;}\newline
\noindent (ii)\textit{\ assuming the Continuum Hypothesis, CH:\newline
}$A\backslash W$\textit{\ is countable for open }$W$\textit{\ with} $\mathbb{%
Q\subseteq }W.$

\bigskip

For set-theoretic background we refer to [BinO8]. We now recall a few of the
classes used to study the adequacy of sets in sustaining (topologically)
`good behaviour' -- notions of adequate size or largeness. These are known
as \textit{gauges}. We compare some of these to help introduce `strange
sets', outside the bounds of the usual standard classification of the size
of a set. Our interest here rests on the following families of subsets of $%
\mathbb{R}$. (Below two subsets are \textit{similar }if they are images
under some (injective) affine function: $f(x)=mx+c$ with $m\neq 0.$)

\bigskip

\noindent \textbf{Definitions. }Put:

\noindent $\mathcal{L}^{+}:=\{A:\lambda (A)>0\}$ with $\lambda $ Lebesgue
measure on $\mathbb{R}$ and $\lambda ^{\ast }$ (below) its outer measure;

\noindent $\mathcal{B}a\mathcal{^{+}}:=\{A:A\text{ is second category and
has the Baire property}\};$

\noindent $\mathcal{SW}:=\{A:d(A)\text{ contains an interval}\};$

\noindent $\mathcal{E}mb:=\{A:\text{ for each finite set $F$, $A$ contains a
set $\widetilde{F}$ similar to $F$}\};$

\noindent $\mathcal{SC}:=\{A:A\text{ is shift-compact}\};$

\noindent $\mathcal{BC}=\{A:\text{ if $f:\mathbb{R\rightarrow R}$ is
additive and bounded above on $A$,}$\newline
then $f$ is linear (i.e. continuous) $\}.$

\bigskip

The last case above is the Ger-Kuczma class $\mathfrak{B}$ of [Kuc,\S 9,10]
(cf. [BinO1]), but with so many families in play our chosen notation above
gives more of a mnemonic (as with $\mathcal{E}mb$ for embedding). The
condition `bounded from above' in $\mathcal{BC}$ can be replaced by `bounded
from below', as $-f$ is additive whenever $f$ is. Note, however, that
replacing either of these by just `bounded', yields a larger class, the
Ger-Kuczma class $\mathfrak{C}$ -- see [Kuc, Th. 9.1.1]. There is in
principle a third Ger-Kuczma class, denoted $\mathfrak{A}$ in [Kuc],
analogous to $\mathfrak{B}$ but referring to mid-point convex functions;
however, it emerges that $\mathfrak{A}=\mathfrak{B}$ [Kuc, Th. 10.2.2].

\bigskip

The first two classes are thoroughly studied in [Oxt] and it is well known
that $\mathcal{L^{+}\cup B}a\mathcal{^{+}}\subseteq \mathcal{SW\cap E}mb$
(for $\mathcal{SW}$ this is the Steinhaus-Weil Theorem, [Oxt, Th. 4.8], cf.
[BinO9]; for $\mathcal{E}mb$ see [Kel] which gives a brief survey, cf. [Sve]
a much earlier survey including the related `Erd\H{o}s similarity problem').
Furthermore, $\mathcal{L^{+}\cup B}a\mathcal{^{+}}\subseteq \mathcal{SC\cap
BC}$ is also well-known (for $\mathcal{BC}$ see $\mathfrak{B}$ in [Kuc, Th.
9.3.3] and for $\mathcal{SC}$ see [BinO5] -- though this goes back to
Kestelman [Kes1] and Borwein and Ditor [BorD]). Note that if $d(A)=A-A$
contains an interval, then $A+A$ need not: see [CrnGH] for an example of a
compact subset $S$ such that $d(S)=S-S$ contains an interval, but $S+S$ has
measure zero.

We recall some recent results and offer some new examples along similar
lines. Firstly,

\bigskip

\noindent \textbf{Theorem MM }([MilM]).

\noindent (i) \textit{There exists a shift-compact set that is concentrated
on }$\mathbb{Q}$\textit{, the rationals}.

\noindent (ii) \textit{There exists a non-measurable shift-compact set.}

\bigskip

These two results appeared recently in [MilM]. We study the possible
topological character of the set in (i) under G\"{o}del's axiom $V=L$ in \S %
5. We continue with some new examples accompanied by earlier results which
they complement. Thus the example in (ii) below is simpler than that in
[Kuc, Th. 9.3.4]. We review the effective nature of the constructions here
in Lemma 2 in \S 5, where we study these examples in the light of $V=L$.
Below $\mathfrak{c}$ denotes the cardinality of the \textit{continuum},
treated here as an initial \textit{ordinal}, as is common in set theory
[Jec], [Kun], cf. [BinO8].

\bigskip

\noindent \textbf{Theorem 3.}

\noindent (i) $\mathcal{SW\not\subseteq BC}$ and $\mathcal{BC\not\subseteq SW%
}$;

\noindent (ii) $\mathcal{E}mb\not\subseteq \mathcal{SC}$, $\mathcal{SC}%
\not\subseteq \mathcal{E}mb$;

\noindent (iii) There exists a set $A\subseteq \lbrack 0,1]$, with $\lambda
^{\ast }(A)=1$ such that $A\not\in \mathcal{E}mb;$

\noindent (iv) There exists a set $A\subseteq \lbrack 0,1]$, $A\cap I$
second category for each closed interval $I\subseteq \lbrack 0,1]$ such that 
$A\notin \mathcal{E}mb;$

\noindent (v) There exists a set $A\subseteq \lbrack 0,1]$, with $\lambda
^{\ast }(A)=1$ such that $A\notin \mathcal{BC}$;

\noindent (vi) There exists a set $A\subseteq \lbrack 0,1]$, $A\cap I$
second category for each closed interval $I\subseteq \lbrack 0,1]$ such that 
$A\notin \mathcal{BC}$;

\noindent (vii) There exists a non-measurable set $A$, $A\in \mathcal{BC}$;

\noindent (viii) There exists $A\in (\mathcal{SW\cap E}mb)\backslash 
\mathcal{SC}$;

\noindent (ix) $\mathcal{SC}\subset \mathcal{SW}$, $\mathcal{SC\subset BC}$.

\bigskip

Before proceeding we mention a beautiful result of Cieselski-Rosenblatt
[CieR,Th. 12] that the Erd\H{o}s and Kakutani [ErdK] set%
\[
C_{EK}:=\{\tsum\nolimits_{k=2}^{\infty }a_{k}/k!:a_{k}\in \{0,1,2,..k-2\}\}, 
\]%
which is a compact perfect set of \textit{measure zero,} is shift-compact.
It was already known [EleS] (cf. [EleT]) that for every perfect set $%
P\subseteq \mathbb{R}$ there is $x\in \mathbb{R}$ with $C_{ES}\cap (x+P)$
uncountable. For further literature on this and related matters see [BarLS].
Notice also that $\mathcal{C}$, the excluded middle-thirds Cantor set in $%
[0,1],$ is compact, and $\lambda (\mathcal{C})=0$, but $\mathcal{C}\in 
\mathcal{SW}$ ($d(\mathcal{C})=[-1,1]$ and $\mathcal{C}+\mathcal{C}=[0,2]$)
and hence $\mathcal{C}\in \mathcal{BC}$. Also $\mathcal{C}\notin \mathcal{E}%
mb$.

\bigskip

\noindent \textbf{Proof of Theorem 3.}

\noindent \textit{Proof of}\textbf{\ }(i). Let $f$ be any \textit{%
discontinuous} additive function on $\mathbb{R}$ (for examples see e.g.
[Kuc, \S 5.2]). Put $A=\{x:f(x)\leq 0\}$; then, as $f$ is bounded from above
on $A$ but not continuous (linear), $A\notin \mathcal{BC}$. However, since $%
f $ is additive and $0\in A$ it is immediate that $d(A)=\mathbb{R}$. (If $%
f(x)>0,$ then $x=0-(-x)\in d(A)).$ So $A\in \mathcal{SW}$. Therefore $%
\mathcal{SW}\not\subseteq \mathcal{BC}$.

For the second part, take a Hamel basis $H=\{h_{\alpha }\}_{\alpha <%
\mathfrak{c}}$, where $\mathfrak{c}$ the cardinality of the continuum as
above, and let 
\[
A:=\{qh_{\alpha }:q\in \mathbb{Q},q\neq 0,\alpha <\mathfrak{c}\}. 
\]%
Then $d(A)$ consists of all $q_{1}h_{\alpha _{1}}-q_{2}h_{\alpha _{2}}$ with 
$q_{1},q_{2}\neq 0$. Fix three distinct $\alpha _{1}$, $\alpha _{2}$, $%
\alpha _{3}$. Then 
\[
q_{1}h_{\alpha _{1}}+q_{2}h_{\alpha _{2}}+q_{3}h_{\alpha _{3}}\notin d(A) 
\]%
whenever $q_{1}\neq 0$, $q_{2}\neq 0$, $q_{3}\neq 0$, and these numbers are
dense in $\mathbb{R}$. Hence $d(A)$ contains no interval. However, $A\in 
\mathcal{BC}$: if $f$ is additive and bounded above on $A$, then $%
f(h_{\alpha })=0$ for every $\alpha ,$ so $f=0,$ and so is vacuously linear
(continuous). $\square $ (i)

\noindent \textit{Proof of} (ii). This falls into two parts.

\noindent \textit{Part} $1.$

We will construct a set $B\in \mathcal{E}mb\setminus \mathcal{SC}$ by
transfinite induction of length $\mathfrak{c}$. Let $F_{\alpha }$, $\alpha <%
\mathfrak{c}$ denote all finite sets of real numbers. Set $B_{0}=F_{0}$. Let 
$\displaystyle{A_{1}:=\{b\pm \frac{1}{n}:b\in B_{0},n\in \mathbb{N}\}.}$
Clearly there exists $\widetilde{F_{1}}$ similar to $F_{1}$, such that $%
\widetilde{F_{1}}\cap (A_{1}\cup B_{0})=\emptyset $. Set $B_{1}=\widetilde{%
F_{1}}\cup B_{0}$.

Now for some $\alpha <\mathfrak{c}$, suppose we have constructed $\langle
B_{\beta },\beta <\alpha \rangle ,$ so that, for each $\beta <\alpha ,$ $%
B_{\beta }=\widetilde{F_{\beta }}\cup \bigcup_{\mathfrak{\gamma }<\beta }B_{%
\mathfrak{\gamma }}$, with $\widetilde{F_{\beta }}$ similar to $F_{\beta }$
and with $\displaystyle{\widetilde{F_{\beta }}\cap (A_{\beta }\cup \bigcup_{%
\mathfrak{\gamma }<\beta }B_{\mathfrak{\gamma }})=\emptyset ,}$ where 
\[
A{_{\beta }=\{b\pm \frac{1}{n}:b\in \bigcup_{\mathfrak{\gamma }<\beta }B_{%
\mathfrak{\gamma }},n\in \mathbb{N}\},} 
\]%
and each $B_{\beta }$ (from construction) has cardinality less than or equal
to that of $\beta $ (when $\beta $ is infinite). Let 
\[
{A_{\alpha }=\{b\pm \frac{1}{n}:b\in \bigcup_{\beta <\alpha }B_{\beta },n\in 
\mathbb{N}\}}. 
\]%
From the inductive hypothesis, $\displaystyle S_{\alpha }:={A_{\alpha }\cup
\bigcup_{\beta <\alpha }B_{\beta }}$ has cardinality less than or equal to
that of $\alpha $, and thus less than $\mathfrak{c}$. So there exists $%
\widetilde{F_{\alpha }}$ similar to $F_{\alpha }$ such that $\displaystyle{%
\widetilde{F_{\alpha }}\cap S}_{\alpha }{=\emptyset }$. (Consider the
similarity $f(t)=at$ with $a\notin S_{\alpha }f^{-1}$ for $f\in F_{\alpha }.$%
) Now define $\displaystyle{B_{\alpha }=\widetilde{F_{\alpha }}\cup
\bigcup_{\beta <\alpha }B_{\beta }}$. If we set $B:=\displaystyle{%
\bigcup_{\alpha <\mathfrak{c}}B_{\alpha }}$, it is routine to verify that $%
B\in \mathcal{E}mb\setminus \mathcal{SC}$.

\noindent \textit{Part} $2.$

We will construct a set $B\in \mathcal{SC}\setminus \mathcal{E}mb$ by
transfinite induction, ensuring that it contains no subset similar to $%
\{1,2,3\}.$

Arrange all the null-sequences in a transfinite sequence $\langle
\{x_{n}^{\alpha }\}:\alpha <\mathfrak{c}\rangle $.

Set $B_{0}=\displaystyle{\{b_{0}\}\cup \{x_{n_{k,0}}^{0}:n_{k,0}\in \mathbb{N%
}\},}$ where $b_{0}=0$ and $\{x_{n_{k,0}}^{0}\}_{n_{k,0}\in \mathbb{N}}$ is
a subsequence of $\{x_{n}^{0}\}_{n\in \mathbb{N}}$ so that $B_{0}$ contains
no set similar to the set $\{1,2,3\}$.

Now suppose that for some $\alpha <\mathfrak{c}$ we have chosen $\langle
B_{\beta }:\beta <\alpha \rangle $ to satisfy%
\[
B_{\beta }=\bigcup_{\mathfrak{\gamma }<\beta }B_{\mathfrak{\gamma }}\cup
\{b_{\beta }\}\cup \{x_{n_{k,\beta }}^{\beta }{:n_{k,\beta }\in \mathbb{N}}%
\}, 
\]%
where $\displaystyle{\{x_{n_{k,\beta }}^{\beta }\}}_{n_{k,\beta }\in \mathbb{%
N}}$ is a subsequence of $\{x_{n}^{\beta }\}_{n\in \mathbb{N}}$, with $%
b_{\beta }$ a real number such that $B_{\beta }$ contains no set similar to
the set $\{1,2,3\}$. Clearly $\displaystyle{\bigcup_{\beta <\alpha }B_{\beta
}}$ has less than $\mathfrak{c}$ elements, so it is easy to verify that we
can choose $b_{\alpha }$ and $\displaystyle{\{x_{n_{k,\alpha }}^{\alpha }\}}%
_{n_{k,\alpha }\in \mathbb{N}}$ a subsequence of $\{x_{n}^{\alpha }\}_{n\in 
\mathbb{N}}$, so that 
\[
B_{\alpha }=\bigcup_{\beta <\alpha }B_{\beta }\cup \{b_{\alpha }\}\cup
\{x_{n_{k,\alpha }}^{\alpha }{:n_{k,\alpha }\in \mathbb{N}}\} 
\]%
contains no set similar to $\{1,2,3\}$.

Finally set $\displaystyle{B=\bigcup_{\alpha <\mathfrak{c}}B_{\alpha }}$.
Then $B$ is shift-compact, as $0\in B$ and each null sequence contains a
subsequence in $B,$ and further $B\notin \mathcal{E}mb$. $\square $ (ii)

\noindent \textit{Proof of} (iii). $A\subseteq \lbrack 0,1]$ satisfies $%
\lambda ^{\ast }(A)=1$ iff $A\cap F\neq \emptyset $ for every closed subset $%
F$ of $[0,1]$ of positive measure. Let $\langle F_{\alpha }$ , $\alpha <%
\mathfrak{c\rangle }$ enumerate the closed subsets of $[0,1]$ of positive
measure. By transfinite induction, we can construct $A=\{x_{\alpha }:\alpha <%
\mathfrak{c}\}$ by picking $x_{\alpha }\in F_{\alpha }$ at each step $\alpha
<\mathfrak{c}$ in such a way that $\{x_{\beta }:\beta \leq \alpha \}$
contains no subset similar to $\{1,2,3\}$. This is possible since at each
step $\alpha <\mathfrak{c}$ we have less than $\mathfrak{c}$ excluded values
for the choice of $x_{\alpha }$, and $F_{\alpha }$ has cardinality $%
\mathfrak{c}$. $\square $ (iii)

\noindent \textit{Proof of} (iv). First notice that for $A\subset \lbrack
0,1]$, $(A)$ and $(B)$ below are equivalent: \newline
\noindent $(A)$ $A\cap F\neq \emptyset $, $\forall F\subseteq \lbrack 0,1]$
with $F$ closed and second category. \newline
\noindent $(B)$ $A\cap I$ is second category $\forall I\subseteq \lbrack
0,1] $, with $I$ a closed interval.

Let $\langle F_{\alpha }$ , $\alpha <\mathfrak{c\rangle }$ enumerate the
collection of second-category closed subsets of $[0,1]$. Again, by
transfinite induction, we can construct $A=\{x_{\alpha }:\alpha <\mathfrak{c}%
\}$ by picking $x_{\alpha }\in F_{\alpha }$ at each step $\alpha <\mathfrak{c%
}$ in such a way that $\{x_{\beta }:\beta \leq \alpha \}$ contains no subset
similar to $\{1,2,3\}$. Then for $A$, $(A)$, or equivalently $(B),$ holds
and $A\notin \mathcal{E}mb$. $\square $ (iv)

\noindent \textit{Proof of} (v). Treating ${\mathbb{R}}$ as a vector space
over ${\mathbb{Q}}$ and with $H$ a Hamel basis, as above, take $A:=Lin_{{%
\mathbb{Q}}}(H\setminus \{h_{0}\})$ the vector subspace generated by $%
H\setminus \{h_{0}\}$ of co-dimension $1$. Then the additive function $f$
generated $by$ taking $f(h_{0})=1$ and $f(h)=0$ for $h\in H\setminus
\{h_{0}\}$ is discontinuous, and bounded on $A.$ So $A\notin \mathcal{BC}$.
Also $\lambda ^{\ast }(A)=1.$ $\square $ (v)

\noindent \textit{Proof of }(vi). Suppose $A$ is the same as in the proof of
(v) so that $A\notin \mathcal{BC}$. We will show that $A\cap I$ is of second
category for every closed interval $I\subseteq \lbrack 0,1]$. Let $I$ be
given. Let $\overline{I}=\displaystyle{\frac{I}{2}}$ with the same centre as 
$I$ and half the length. Now 
\[
\overline{I}=\bigcup_{r\in \mathbb{Q}}[(A+rh_{0})\cap \overline{I}%
]=:\bigcup_{r\in \mathbb{Q}}T_{r.} 
\]%
Since $\overline{I}$ is of second category, at least one $T_{r}$ is of
second category, $T_{\overline{r}}$ say. That is, $(A+\overline{r}h_{0})\cap 
\overline{I}$ is of second category.

Take $x\in T_{\overline{r}}$ . Since $A$ is dense in $\mathbb{R}$, $%
\overline{r}h_{0}$ can be written as $\overline{r}h=\overline{a}+\epsilon $
with $\overline{a}\in A$ and $|\epsilon |<\frac{|I|}{5}$, and hence $x=a+%
\overline{a}+\epsilon $ for some $a\in A$. So $T_{\overline{r}}\subseteq
(A+\epsilon )\cap \overline{I},$ and so $(A+\epsilon )\cap \overline{I}$ is
of second category. This implies that $A\cap I$ is of second category, being
a translate by $-\epsilon $ of the set $(A+\epsilon )\cap \overline{I}$,
completing the proof. $\square $ (vi)

\noindent \textit{Proof of} (vii). Take $B\subseteq (0,1)$, with $B$
non-measurable. Then $A=B\cup \lbrack 1,2]$ is automatically non-measurable,
and in $\mathcal{SC}$, and so in $\mathcal{BC}$, by Darboux's theorem (see
e.g. [BinO5]). $\square $ (vii)

\noindent \textit{Proof of} (viii). Let $A$ be the set of Theorem MO above
(constructed in the proof of Theorem 9 in [MilO]). Then $A\in \mathcal{%
SW\cap E}mb$. We will show $A$ is not shift-compact. Suppose otherwise, and
consider the null sequence $\displaystyle{(-1/n)}$. We show that 
\[
A\cap \bigcap_{k=1}^{\infty }(A-\frac{1}{n_{k}})=\emptyset , 
\]%
for \textit{every} subsequence $\displaystyle{(-1/n}_{k}{)}$, so
contradicting that $A\ $is shift-compact. So suppose the intersection above
is non-empty for some subsequence $\displaystyle{(-1/n}_{k}{)}$. Then, as $A$
is assumed shift-compact, there exist $a\in A$ such that $\displaystyle{%
a+(1/n}_{k})\in A$ for all $k$, which is impossible by Th. MO(i) (i.e. (c)
in Theorem 9 of [MilO]). $\square $ (viii)

\noindent \textit{Proof of} (ix). It is a corollary of earlier parts,
already proved, that these two inclusions are proper: thus the first being
proper follows from (viii). Both $\subseteq $-inclusions are well-known: see
[BinO1, Th.1] for the first and [BinO5] for the second. For completeness, we
recall the inclusion proofs here, as they are short (and needed together
below).

Suppose $A\in \mathcal{SC}$. We claim that $[0,\delta )\subseteq d(A)$ for
some $\delta >0$. Otherwise, there exists a null sequence $y=\{y_{n}\}_{n\in 
\mathbb{N}}$, with $y_{n}\notin d(A)$ for each $n\in \mathbb{N}$. Since $%
A\in \mathcal{SC}$, there exists a subsequence $\{y_{n_{k}}\}_{n_{k}\in 
\mathbb{N}}$, and $a\in A$ such that $a+y_{n_{k}}\in A$ for all $k$. Hence $%
y_{n_{k}}=(a+y_{n_{k}})-a\in d(A)$ for all $k$, a contradiction. Thus $%
\mathcal{SC}\subseteq \mathcal{SW}$.

Now we show $\mathcal{SC}\subseteq \mathcal{BC}$. Suppose otherwise. Then
there is a shift-compact set $A\notin \mathcal{BC}$ and an additive function 
$f$ on $R$ that is discontinuous but bounded above on $A$. By Darboux's
theorem, there exists a \textit{null} sequence $\{y_{n}\}_{n\in \mathbb{N}}$
with $f(y_{n})\longrightarrow \infty $ (as otherwise $f$ is locally bounded
at $0$, and so continuous again by Darboux's theorem$).$ Since $A$ is
shift-compact, there exists $a\in A$ and a subsequence $\{y_{n_{k}}\}_{n_{k}%
\in \mathbb{N}}$ such that $a+y_{n_{k}}\in A$ for all $k\in \mathbb{N}$.
Then $f(a+y_{n_{k}})=f(a)+f(y_{n_{k}})\longrightarrow \infty $, a
contradiction since $f$ is bounded from above on $A$.

By the earlier part of this proof, $\mathcal{SC}\subseteq \mathcal{SW\cap BC}
$; but, as in (i), $\mathcal{BC\not\subseteq SW}$, so again this is a proper
inclusion. $\square $ (ix)

\bigskip

An alternative example for (viii) is provided by [MilM], cf. Th. MM (ii)
above. See also \S\ 6.2 below on Sierpi\'{n}ski sets. We stress that the
inclusions mentioned in the Theorem 3 above are all proper, as shown in the
proofs.

\section{Singular sets of good character}

In this section we reconsider an earlier counter-example and show that under 
$V=L$ it will have good character: it will be co-analytic.

We recall from \S 1 that a subset $T$ of the reals is \textit{shift-compact}
if for any null sequence $z_{n}\rightarrow 0$ there is $t\in T$ such that $%
t+z_{m}\in T$ for infinitely many $m.$ We refer to $t$ as a `translator into 
$T$ for $z$'.

Recall that a set $S$ is \textit{concentrated} on the rationals $\mathbb{Q}$
if it is uncountable and for every open set $W\supseteq \mathbb{Q}$ the set $%
S\backslash W$ is countable [Rog1, \S 2.3]. Such a set is of strong measure
zero. Under the assumption that less than $\mathfrak{c}$ meagre sets have
meagre union, the first two authors have shown in [MilM] (cf. Th. MM in \S 4
above) that there is a set concentrated on $\mathbb{Q}$ which is
shift-compact. To discuss a refinement of this result involving effective
aspects, we recall the (effective) \textit{analytical hierarchy} of
predicates in the language of set theory (i.e. with the non-logical symbol $%
\in ,$ cf. \S 3.1) concerned with numbers (members of $\omega )$ and `reals'
(represented by number sequence in $\omega ^{\omega }$). See e.g. [Rog2,
Part 4] or [BinO8, \S 8 The syntax of analysis] for background. Write these
with all quantifiers $\forall x$ and $\exists x$ (ranging over reals $x$) at
the front, followed by an arithmetical predicate (this can be done assuming
the Axiom of Dependent Choices, DC); then list and name the predicates
according to the starting quantifier and the number alternations (between $%
\forall $ and $\exists ,$ universal and existential) binding all the
variables. Thus, as above in \S 3.1, a (lightface) $\Sigma _{1}^{1}$
predicate has just one existential and $\Pi _{1}^{1}$ has just one
universal; $\Sigma _{2}^{1}$ has the $\exists \forall $ format, etc. If a
free variable \textit{parameter} $x\in (0,1)$ is allowed in the predicate
(with $x,$ regarded via its binary expansion as a function with domain $%
\mathbb{N}$, not necessarily effectively defined), this is recognized by
bold-face lettering, yielding a hierarchy that is `relativized' in the
parameter (permitting relative effectiveness [Kec, V.40B]). Here $\mathbf{%
\Sigma }_{1}^{1}$ corresponds to classical analytic sets and $\mathbf{\Pi }%
_{1}^{1}$ to the co-analytic sets: an arbitrary open set in the line can be
coded by a not necessarily recursive sequence of the rational-ended basic
open intervals it contains.

The first two authors' result amends a classical construction of a
concentrated set using transfinite induction -- so that, as first noted by
Kuratowski [Kur] -- under $V=L$ such a set would be $\Delta _{2}^{1}=\Sigma
_{2}^{1}\cap \Pi _{2}^{1}.$ In fact, under $V=L,$ as [ErdKM, Th. 13] have
shown, with careful monitoring of the effectiveness of constructions, a set $%
S$ concentrated on $\mathbb{Q}$ can be constructed which is $\Pi _{1}^{1}.$
The underlying reason for the character improvement is that their
construction is based on combinatorial analysis that is suitably `effective'.

We will similarly demonstrate an effective construction of a translator for
a null sequence $z$ into any dense $\mathcal{G}_{\delta }$ set $T$, when $%
T=\bigcap\nolimits_{n}G_{n}$ with each $G_{n}$ open and containing $\mathbb{Q%
}.$ This uses an effective enumeration of $\mathbb{Q}$ and the fairly recent
constructive proof of shift-compactness [MilO]. We regard this as a
geometric counterpart to the more combinatorial argument of [ErdKM],
establishing the following

\bigskip

\textbf{Theorem 4.} \textit{Under }$V=L,$\textit{\ there is a }$\mathbf{\Pi }%
_{1}^{1}$\textit{\ subset of the reals which is concentrated on} $\mathbb{Q}$
\textit{and is shift-compact. }

\bigskip

The result is not altogether surprising. In the language of Turing
reducibility (below), Vidny\'{a}nszky [Vid] captures the general procedure
of adapting a construction of a set $S$ in a Polish space by transfinite
induction under the assumption $V=L$ to yield a \textit{coanalytic version} $%
C$ of $S$ in the following formulation, a result implied by $V=L$:

\bigskip

\noindent \textbf{Theorem V} ([Vid, Th. 1.3]). \textit{Assume V=L. For }$B$%
\textit{\ an uncountable Borel subset of an arbitrary Polish space, if} 
\newline
\qquad (i)\textit{\ }$F$\textit{\ is a co-analytic subset of }$M^{\leq
\omega }\times B\times M$\textit{\ with }$M\in \{R^{n},2^{\omega },\wp
(\omega ),\omega ^{\omega }\},$\textit{\ and}\newline
\qquad (ii)\textit{\ for all }$A\in M^{\leq \omega }$\textit{, }$p\in B,$%
\textit{\ the vertical section }$F(A,p)\subseteq M$\textit{\ of }$F$\textit{%
\ is (upwards) cofinal in the ordering }$\leq _{T}$\textit{of
Turing-reduciblity}

\noindent \textit{-- then there exists a co-analytic set }$C$\textit{\ that
is `compatible with }$F$\textit{'.}

\bigskip

Here $M^{\leq \omega }$ denotes the countable subsets of $M,$ and we recall
that $x\leq _{T}y$ for $x,y\in M$ (read: `$x$ is Turing reducible to $y$'),
if $x$ can be effectively computed from $y$ (more exactly: there exists a
Turing machine which computes $x$ from the input $y$).

Rather than apply Th. V, which shadows [ErdKM], we have ourselves shadowed
[ErdKM] in the preamble to the proof of Theorem 4 in an exposition of the
tools from logic -- which we hope analysts will find congenial -- thereby
clarifying the nub of the result. We rely on specified background from an
analyst-friendly source: [BinO8].

The proof of Theorem 4 is given below, as indicated. We preface that now
with a discussion of its salient features, in particular on its reliance on
Kleene's theorem below (which gives a circumstance when an existential
quantifier can be converted to a universal one).

\bigskip

\textit{Proof of Theorem 4 preamble: proof strategy.}

We need to refer to the (metamathematical -- `external' to the discourse in
the language) semantic relation $\models $ of satisfaction/truth (below),
due to Tarski (see [Tar1,2], cf. [BelS, Ch. 3 \S 2], cf. [BinO8]), which is
read as `models', or informally as `thinks' (adopting a common enough
anthropomorphic stance). A formula $\varphi $ of $LST$ with free variables $%
x,y,...,z$ may be interpreted in the structure $\mathcal{M}:=\langle M,\in
_{M}\rangle $ (with $\in _{M}$ now a binary set relation on the set $M$) for
a given assignment $a,b,...,c$ in $M$ for these free variables, and one
writes%
\[
\mathcal{M}\models \varphi (x,y,...,z)[a,b,...,c],\text{ or by abbreviation }%
M\models \varphi \lbrack a,b,...,c], 
\]%
if the property holds; this requires an induction on the syntactic
complexity of the formula starting with the atomic formulas (for instance,
the atomic case $x\in y$ is interpreted under the assignment $a,b$ as
holding iff $a\in _{M}b$).

We recall also that \textit{Skolemization} of a formula of $LST,$ say $\phi (%
\bar{x},\bar{\tau})$ with $\bar{x}$ a finite list of \textit{free variables}
and $\bar{\tau}$ a finite \textit{list of ordinals}, is the elimination of
all quantifiers by\newline
\qquad (i) replacing existential quantifiers with functions appointing a
`witness' of an asserted existence (from among the available instances,
assuming any exist), and\newline
\qquad (ii) making free the variables previously bound by universal
quantifiers (for which see [Hod, Ch. 3, p. 71], cf. [ManW, p. 87]).

This process yields an `equi-satisfiable' (equivalent under $\models $)
quantifier-free formula $\Phi (\bar{x},\bar{\tau},\bar{y},\bar{f}),$
involving a further finite list of free variables $\bar{y}$ and finite list
of function symbols $\bar{f}$ (the Skolem functions for $\phi )$ arising
from the Skolemization, such that%
\[
\exists \bar{f}\forall \bar{y}[\phi \rightarrow \Phi \wedge \psi ] 
\]%
is a theorem of predicate logic (suppressing here the various lists $\bar{x},%
\bar{\tau},\bar{y},\bar{f}$); here $\psi $ is a certain (known) sentence
such that if $M$ is a transitive set and $\psi $ holds in $M,$ then $M$ is
an $L_{\alpha }$.

The structure $\langle L_{\alpha },\in \rangle $ can be equipped with 
\textit{canonical Skolem functions} through always appointing `witnesses' as
above that are earliest under the well-ordering $<_{L}$ of \S 3.1 above. Say
that $L_{\alpha }$ is \textit{point-definable} if its Skolem hull (smallest
set including $L_{\alpha }$ and closed under the iteration of all its
canonical Skolem functions) is isomorphic to $L_{\alpha }.$ (Such $L_{\alpha
}$ exist for unboundedly many $\alpha $ in $\omega _{1}$ -- for proof see
[EngMS, Proof of Th. 2.6].) Performing the canonical Skolemization of $%
L_{\alpha }$, one may define a relation $E_{\omega }$ on $\omega ,$
recursive in the set of all (first-order) sentences true in $L_{\alpha }$
(known as the `theory of $L_{\alpha }$', denoted $Th(L_{\alpha })$) such
that $(\omega ,E_{\omega })\approx (L_{\alpha },\in ),$ where $\approx $
denotes isomorphism (see [ManW, p. 87]). Bearing in mind its definability, $%
E_{\omega }\in L_{\alpha +3},$ since $Th(L_{\alpha })\in L_{\alpha +2}.$

Consequent on the effective combinatorics used in the transfinite inductive
construction in [ErdKM], membership of the singular set $S$ constructed
there can be expressed by a formula, denoted $\mathcal{S}(.)$ (with one free
variable), in such a way that if $x=x_{\alpha }\in L$ is selected
inductively by reference to a point-definable $(L_{\alpha },\in )$ and to
the ordering $<_{L},$ then one constructs, recursively in $x$ and in $%
Th(L_{\alpha }),$ a countable set $M$ and a relation $E_{M}$ on $M$ such
that $(M,E_{M})\models \mathcal{S}(x)$ (i.e., the sentence $\mathcal{S}(x)$
holds in the structure $M$). Taking $z$ to code $Th(L_{\alpha }),$ a real $%
\mu $ may be constructed from $z$ to code the set $M$ and the relation $%
E_{M} $ on $M;$ when done effectively the real $\mu $ is called \textit{%
recursive in} $z.$ Indeed, $(M,E_{M})$ may be constructed to be isomorphic
to $(L_{\alpha +\omega },\in ),$ cf. [EngMS, Th. 2.6, p. 209].

To verify the $\Pi _{1}^{1}$ character of the set $S,$ [ErdKM] relies on 
\textit{Kleene's theorem} from recursion theory (for which see e.g. [Sac,
Lemma 3.1.III] and the formal proof below) that the existential quantifier
over the `reals recursive in $z$' (and, more generally, to reals in the set $%
HYP(z)$ that are `hyperarithmetic in $z$') may in fact be rendered as a
universal quantifier ranging over all the reals. (See [Sac, Lemma 3.1.III],
or [ManW, 4.19]; note that there are countably many reals hyperarithmetic in 
$z.$) Now the satisfaction relation `$M\models \mathcal{S}(x)$' when applied
to countable models $M$ is $\Delta _{1}^{1}$ as a predicate involving the
real number $\mu $ coding $M$, as above (see e.g. [ManW,1.20]), so it is in
particular $\Pi _{1}^{1}$. Now $x\in S$ iff 
\[
\exists \mu \in HYP(x)[\mu \approx (L_{\alpha +\omega },\in )\&\mu \models 
\mathcal{S}(x)], 
\]%
which is $\Pi _{1}^{1}$ in $x$ (so $\mathbf{\Pi }_{1}^{1}$) by Kleene's
theorem. The main task in the formal proof Theorem 4 below is analogous: to
convert the informal statement \textquotedblleft $x\in X$\textquotedblright\
to a formula $\mathcal{S}(x),$ the idea being to recover it, as in [ErdKM]
above, from the (somewhat circuitous) definition:%
\[
x\in X\iff \exists M\in HYP(x)[M\approx L_{\alpha +\omega }\&M\models (\text{%
\textquotedblleft }x\in X\text{\textquotedblright })]. 
\]%
Once this is done, one may ostensibly again apply Kleene's theorem, but
needs to check that the satisfaction clause (the last clause in the display
above) does not degrade the descriptive character of the entire contents of
the square brackets. One needs the final clause to be $\Pi _{1}^{1}.$
However, the satisfaction relation $M\models P(x)$ arising here is defined
(by induction on the complexity of the predicate) only for predicates $P(x)$
written in $LST$ subject to the restriction that constants involved in $P(x)$
(including $x$ itself) may name only elements of $M.$ (This ensures that
these constants have interpretations in $M;$ in particular, $M$ needs to
contain $x$).

\bigskip

\textit{Proof of Theorem 4: Formal proof.} Assuming $V=L$, we have $%
I\subseteq L_{\omega _{1}.}$ For $\alpha <\omega _{1},$ let $L_{\alpha }$ be 
\textit{point-definable} (as in the preamble above). Select a dense subset $%
D\subseteq I$ such that some $d\in I$ is its recursive enumeration $%
d=\{d_{n}\}_{n\in \mathbb{N}},$ with $d_{n}$ the $n^{\text{th}}$ projection
of $d.$ Put $\mathbb{G}:=\{x\in I:G(x)\supseteq D\}.$ For $x\in I\cap
L_{\alpha }$ with $G(x)$ containing $D,$ the set $I\backslash G(x)$ is
nowhere dense and so 
\[
M_{\alpha }:=\bigcup\nolimits_{x\in L_{\alpha }\cap \mathbb{G}}(I\backslash
G(x)) 
\]%
is meagre, as $L_{\alpha }$ is countable. Put $B_{\alpha }:=I\backslash
M_{\alpha }.$ As there are countably many null sequences in $L_{\alpha },$
there is $t\in I\backslash L_{\alpha }$ such that:\newline
\qquad (a) $t_{n}\in B_{\alpha }$ for each $n,$ and\newline
\qquad (b) for each null sequence $z\in I\cap L_{\alpha }$ there is $m=m(z)$
and $N=N(z)\in \mathbb{N}$ such that $t_{m}+z_{n}\in B_{\alpha }$ for $%
n>N(z).$ As above such a $t$ lies in $L_{\alpha +\omega }.$

Proceed as in [ErdKM], and define $X$ to be the set of all $x\in I$ such
that there exist:\newline
\qquad (i) a limit ordinal $\alpha $ such that $L_{\alpha }$ is
point-definable and satisfies \textrm{ZF}$^{-}$,\newline
\qquad (ii) $E\subseteq \omega \times \omega $ recursive in $x$ such that $%
(\omega ,E)$ is isomorphic with $(L_{\alpha },\in ),$\newline
\qquad (iii) $x$ is the first element of $I\backslash L_{\alpha }$
satisfying (i) and (ii) and (a) and (b) above.

As in the preamble we are to apply \textit{Kleene's theorem} that, for
arithmetic $A(n,f)$ (here $n$ ranges over $\mathbb{N}$ and $f$ over $\mathbb{%
N}^{\mathbb{N}}$), the predicate $(\exists f\in HYP)A(n,f)$ is $\Pi _{1}^{1}$
-- see e.g. [Sac, Lemma 3.1.III]], (In fact, the Spector-Gandy Theorem [Sac,
Th. 3.5] asserts that this format characterizes $\Pi _{1}^{1}.)$ We need to
verify that the defining clauses (a) and (b) and (i)-(iii) are satisfied in
the model $(\omega ,E)$.

To this end, we note that when the satisfaction relation $\models $ is
restricted to a $\Sigma _{1}^{1}$ predicate $P(x)$, it is a $\Pi _{1}^{1}$
relation (in $x)$ -- see [Sac, Lemma 4.5.III]. Alternatively, for the
relation to be $\Pi _{1}^{1}$ the predicate $P(x)$ needs to be a \textit{%
ranked} one, i.e. an ordinal bound $\alpha <\omega _{1}^{L}$ must be placed
on the ranges of the analytic quantifiers and on the free variables
appearing in $P(x)$. (Here $\omega _{1}^{L}$ denotes the ordinal recognized
in $L$ as the first uncountable; it is in fact countable -- see e.g. [Dra, 
\S 8.4] or [BinO8, \S 5.2].)

With this in mind, we check that the defining clauses (a) and (b) and
(i)-(iii) are ranked. Conditions (i) and (ii) are manifestly ranked, as will
be (iii) provided also (a) and (b) are. For (a) one has 
\[
y\in B_{\alpha }\iff y\notin M_{\alpha }\iff (\forall z\in L_{\alpha }\cap
I)[z\in \mathbb{G\rightarrow }\text{ }y\in G(z)], 
\]%
and whilst this is $\Pi _{1}^{1}$ (rather than $\Sigma _{1}^{1})$ the
quantifier is bounded to $L_{\alpha };$ so this is actually of ambiguous
class $\Delta _{1}^{1},$ i.e. both $\Pi _{1}^{1}$ and $\Sigma _{1}^{1}$ in
the codes (`notations') for $\alpha $. For (b) note that%
\[
(\forall z\in Null\cap L_{\alpha })\exists m\exists k\forall
n>k[x_{m}-z_{n}\notin M_{\alpha }], 
\]%
where $Null$ stands for the set of null sequences (see the Appendix), and
this is again $\Pi _{1}^{1}$, but nevertheless the quantifier is bounded to $%
L_{\alpha },$ so is again $\Delta _{1}^{1}$ in the codes for $\alpha .$ $%
\square $

\bigskip

Theorem 4 above offers a co-analytic version (under $V=L)$ of the example of
Th. MM(i), but not of (ii), as co-analytic sets are measurable. Co-analytic
versions may likewise be obtained for the examples of Theorem 3 above (again
except for the non-measurable example of (vii)). This is a consequence of
the effective nature of the constructions used:

\bigskip

\noindent \textbf{Lemma 2.} (a)\textit{\ Take} $S\subseteq \mathbb{R}$ 
\textit{countable; then all but at most countably many affine
transformations }$f(t):=at+b$\textit{\ map any finite set }$F$\textit{\ to
the complement of }$S,$\textit{\ and in particular:}

\noindent (i) \textit{with }$b\notin B:=S+S-S$\textit{\ countable, }$%
f(\{1,2,3\})\nsubseteqq S;$

\noindent (ii)\textit{\ if further }$a\notin T\cap (T/2)\cap (T/3)$\textit{\
with }$T:=S-B$\textit{\ countable, then} $f(\{1,2,3\})\subseteq \mathbb{R}%
\backslash S.$

\noindent (b) \textit{For non-meagre closed }$F$\textit{\ there is a closed
nwd set }$N$\textit{\ with} $F=N\cup G(\{n:I_{n}\subseteq F\}).$

\bigskip

\noindent \textbf{Proof} (a): We consider the case $F=\{1,2,3\},$ as typical
since the generalization is just a tedious exercise in linear algebra. As $%
f(\{1,2,3\})=\{a+b,2a+b,3a+b\},$ if (i) fails, then $b:=f(1)+f(2)-f(3)\in
S+S-S$; from here (ii) is immediate as $f(i)\in S$ iff $a\in (S-B)/i.$ $%
\square _{\text{(a)}}$

\noindent (b) $G:=\bigcup \{I_{n}:I_{n}\subseteq F\}$ is the interior of $F$
and so $F\backslash G$ is closed and nwd. (This decomposition also appears
in [BinO11, Th. 6M(b)] $\square _{\text{(b)}}.$

\bigskip

\noindent \textbf{Theorem 4}$^{\prime }$\textbf{.} \textit{Under }$V=L,$%
\textit{\ there is a }$\mathbf{\Pi }_{1}^{1}$\textit{\ subset of the reals
that is shift-compact but not in }$\mathcal{E}mb$\textit{, and likewise
there is a }$\mathbf{\Pi }_{1}^{1}$\ \textit{subset in }$\mathcal{E}mb$ 
\textit{which is not shift-compact. Indeed, under }$V=L,$\textit{\ all the
examples in the proof of Theorem 3, save for 3}(vii), \textit{have }$\mathbf{%
\Pi }_{1}^{1}$\textit{\ versions.}

\bigskip

\noindent \textbf{Proof of Theorem 4}$^{\prime }.$

\noindent (i) Here by [MilA2] (cf. [Vid]) there is a co-analytic Hamel basis 
$H,$ and so the set $A=\bigcup\nolimits_{q\in \mathbb{Q}}qH$ is also
co-analytic and in $\mathcal{BC}$, but not shift-compact, since it fails to
have the Steinhaus-Weil property.

\noindent (ii) Part 1. The finite subsets of $\mathbb{R}$ are in an
effective 1-1 correspondence with $\mathbb{R}$ and effective choices of
affine similarities may be made on the basis of Lemma 2(a).

\noindent Part 2. This follows from Lemma 2(a).

\noindent (iii) Each $F_{\alpha }$ may be coded, as in \S 3.1, by its
complement $G(a):=[0,1]\backslash F_{\alpha }$ with $|G(a)|<1.$ The latter
property is arithmetical, being equivalent to the existence of a rational $%
q<1$ with $|\bigcup\nolimits_{m\in F}G(a(m))|<q,$ for all finite $F\subseteq 
\mathbb{N}$.

\noindent (iv) By Lemma 2(b), each set $F_{\alpha }$ may be expressed in the
form $G_{\alpha }\cup N_{\alpha }$ with $N_{\alpha }$ closed nwd and $%
G_{\alpha }$ coded as $G(a)$ where $a(n)=1$ iff $I_{n}\subseteq F_{\alpha }.$
We may thus pick $x\in G(a)$ avoiding both $N_{\alpha }$ (as in Th. 1E
above) and the further countable set generated by the choices made earlier
in the transfinite induction.

\noindent (v) $A$ is co-analytic as in (i).

\noindent (vi) This refers to the same set as in (iii).

\noindent (vii) Co-analytic sets are measurable [Rog2, Th. 2.9.2] [Kec,
29.7].

\noindent (viii) This is covered by Th. 4.

\noindent (ix) This refers back to (i). $\square $

\section{Complements: related singular sets}

\noindent 1. \textit{Luzin sets. }Recall that a Luzin set $L$ is an
uncountable set which meets every meagre set in at most a countable set. A
Luzin set does not have BP.

If $L$ were Baire it would be non-meagre as $L$ is uncountable. But then $L\ 
$is co-meagre on some rational-ended interval $I,$ so w.l.o.g. is a dense $%
\mathcal{G}_{\delta }$ on $I,$ and so contains an uncountable meagre set, a
contradiction.

Hence $L$ cannot be analytic or co-analytic. This means that the $V=L$
construction of \S 5 above cannot be improved to yield a Luzin set.

Marczewski observed in 1938 that a set $L$ is Luzin iff $L$ is uncountable
and is concentrated on every countable dense set. (Clear, since $L$ is Luzin
iff for every dense open $G,$ $L\backslash G$ is countable.) As such, $L$ is
of strong measure zero (SMZ). (Being of measure zero, it is Lebesgue
measurable.)

The first two authors' example in Th. MM(i) of \S\ 4 may be made Luzin, so
despite being SMZ it is shift-compact.

\noindent 2. \textit{Sierpi\'{n}ski sets.} Recall that a Sierpi\'{n}ski set
is an uncountable set which meets every measure zero set in at most a
countable set. It is known from work of (Szpilrajn-)Marczewski and
Kuratowski that a Sierpi\'{n}ski set $S$ is not only meagre, but in fact
perfectly meagre (i.e. $S\cap P$ is meagre in $P$ for any perfect set $P$)
-- see e.g. A. Miller's survey article for a proof [MilA1, Th. 4.1 and 5.2].

If $S$ were measurable, then it would be of positive measure, as $S$ is
uncountable. So $S\ $then contains a compact subset of positive measure,
inside which there exists an uncountable set of measure zero -- just repeat
the construction of the Cantor set. This contradicts the defining property
of $S.$ So $S$ is not measurable.

In fact its complement, $\mathbb{R}\backslash S,$ also non-measurable, is
shift-compact by virtue of being co-meagre. As this is thematic, we give a
direct proof based on KBD of the following.

\bigskip

\noindent \textbf{Proposition 3.} \textit{If }$S$\textit{\ is a Sierpi\'{n}%
ski set and }$z_{n}\rightarrow 0$\textit{\ is null, then for quasi all }$t$%
\textit{\ one has} $t+z_{n}\subseteq \mathbb{R}\backslash S$ \textit{for all 
}$n.$

\bigskip

\noindent \textbf{Proof.} Choose $H$ a dense $\mathcal{G}_{\delta }$ of zero
measure containing the points $z_{n}.$ As $S$ is a Sierpi\'{n}ski set, $%
D:=S\cap H$ is countable, and $S\subseteq D\cup (\mathbb{R}\backslash H)$.
Now $T:=(\mathbb{R}\backslash H)\cup D$ is meagre, so $(\mathbb{R}\backslash
T)=H\backslash D$ is co-meagre. By KBD,\ for quasi all $t\in H\backslash D$
one has $t+z_{n}\subseteq H\backslash D\subseteq \mathbb{R}\backslash S$ for
all $n.$ $\square $

\bigskip

The result above also follows from a stronger result of Jasi\'{n}ski and
Weiss [JasW], concerning shifting a null $\mathcal{F}_{\sigma }$ (`measure
zero' null) rather than a null sequence, and from Carlson [Car], who also
studies associated $\sigma $-ideals. See \S 6.7 below.

\noindent 3. \textit{Characterization of Strong Measure Zero sets (SMZ). }

A set $X$ is of strong measure zero if for each sequence $\{\delta
_{n}\}_{n\in \mathbb{N}}$ with each $\delta _{n}>0$ there is a corresponding
squence of intervals $\{I_{n}\}_{n\in \mathbb{N}}$ covering $X$ with each $%
I_{n}$ of length at most $\delta _{n}.$ Such sets $X$ are characterized by
the property that, for each \textit{meagre} set $H,$ there is $x$ with $%
X\cap (x+H)=\emptyset .$ See [MilA1, Th. 3.5].

Carlson [Car, Th. 2.1] shows that under MA$_{\kappa }$ (Martin's Axiom at $%
\kappa <\mathfrak{c}$) these sets are closed under unions of size $\kappa .$
In particular, this is so for countable unions. He also shows that no
perfect set can be covered by such a countable union (in both $\mathbb{R}$
and in the Cantor space).

In this context we recall the contrasting property [EleS] (cf. [EleT]) of
the Erd\H{o}s-Kakutani set $C_{ES}$ of \S\ 4 that, for every \textit{perfect}
set $P,$ there is $x\in \mathbb{R}$ with $C_{ES}\cap (x+P)$ uncountable.
Compare also \S \S\ 6.4 and 6.7 below.

For further characterizations of \textit{SMZ} see [GalMS]. We mention one
such which is thematic for the present context. Here the target sets $T$ for
embeddings are dense $\mathcal{G}_{\delta }$-sets. Embeddings which are
performed simultaneously in any neighbourhood by a perfect subset of any
such $T$ of a fixed set $Z$ into $T$ characterize those sets $Z$ that are
strongly measure zero. Since any countable set is strongly of measure zero
this result includes `simultaneous embeddings' of a null sequence.

\noindent 4. \textit{Strongly meagre (strong first category). }By analogy
with \textit{SMZ}, a set $X\ $is \textit{strong first category }if for any
measure-zero set $N$ there is $t$ with $X\cap (t+N)=\emptyset .$ See [BarS].

\noindent 5. \textit{Consistency results. }Laver has proved in [Lav1,2] that
it is consistent that every strong measure zero set is countable. Carlson
[Car] shows that likewise it is consistent that every \textit{strong first
category} set is countable.

\noindent 6. \textit{Luzin/Sierpi\'{n}ski sets versus SMZ. }Every Luzin set
has strong measure zero -- see 1 above (this is (Szpilrajn-)Marczewski's
observation).\textit{\ }Bartoszy\'{n}ski and Judah [BartJ, Th. 2] show that,
under the continuum hypothesis CH, every Sierpi\'{n}ski set is a union of at
most two \textit{SMZ}\ sets.

\noindent 7. \textit{Carlson's }$\sigma $\textit{-ideals. }Extending the 
\textit{SMZ} idea, Carlson [Car, Th. 5.7] proves that each of the following
families of sets forms a $\sigma $-ideal:

(i) those sets $X$ with the property that for every a \textit{meagre} set $M$
there is $t$ such that $X\cap (t+M)=\emptyset $;

(ii) those sets $X$ with the property that for every a\textit{\ null} $%
\mathcal{F}_{\sigma }$ set $H\ $there is $t$ such that $X\cap
(t+H)=\emptyset $. (Equivalently, for every $\mathcal{G}_{\delta }$ set $G$
of full measure (=co-null), $X$ is covered by some translate of $G.$)

\bigskip

\noindent 8. \textit{Effective versions and coding. }The proof of Theorem 4
above relied on the ability to refer to various subsets of the real line,
especially open sets, in terms of `codes'. Our canonical sources there were
[Kec, Ch.V] on the analytical hierarchy (and the note [Kec, V.40B] on
classical versus effective descriptive set theory), and our recent survey
[BinO8], and for coding the wide-ranging use in [Sol, II.1.1, 25-33] and the
much more minimal amount in [FenN, \S\ 2, p. 93].

\bigskip

\textbf{Apendix.}We begin with some notation.

Let $\{I_{n}\}_{n\in \mathbb{N}}$ enumerate (constructively) all the
rational-ended intervals, with $I_{n}=(l_{n},r_{n})$. Write $\mathbb{M}$ for
the odd natural numbers; for $a\subseteq \mathbb{N}$ we may extract an $n^{%
\text{th}}$ canonical subset of $a$ and also an open set naturally `coded'
by $a$ by setting:%
\[
a(n)=a\cap \{2^{n}m:m\in \mathbb{M}\},\qquad G(a):=\bigcup\nolimits_{n\in
a}I_{n}. 
\]%
We identify $a\subseteq \mathbb{N}$ with the real number in $\{0,1\}^{%
\mathbb{N}}$ whose binary expansion is the indicator function of $a$. Thus $%
\{a:m\in a\}$ is open (being the set of reals with $m$-th binary digit =1).

\bigskip

\noindent \textbf{Examples. }1.\textbf{\ }Say that $z\in I$ represents a 
\textit{null sequence}, briefly Null$(z)$, if for each $k$ there is $n$ so
that $x_{m}|k=0_{k}$ for all $m\geq n$ (so $z_{n}\rightarrow 0).$ Thus%
\[
\text{Null}(z)\leftrightarrow \forall k\exists l(\forall n\geq l)(\forall
m)[|z(2^{n}(2m+1))|<1/k]. 
\]

\noindent 2. Let $D:=\{d_{n}:n\in \mathbb{N}\}$ enumerate effectively a
subset dense in $I.$ By abuse of notation, say that $x$ \textit{contains} $D$
when $G(x)\supseteq D,$ i.e. for each $n$ there is $m$ with $d_{n}\in
\varphi (x(m)).$ We denote the set of such $x$ by $\mathbb{G}$. Since%
\[
x\in \mathbb{G}\iff (\forall n\in \mathbb{N)(}\exists m\in \mathbb{N)(}%
\exists k\in \mathbb{N)}[d_{n}\in \varphi (k)\text{ and }k=x(m)], 
\]%
this is an arithmetic relation which is (light-faced) $\mathbf{\Pi }%
_{2}^{0}. $

\bigskip

\section{\textbf{References }}

\noindent \lbrack Anc] F.D. Ancel, An alternative proof and applications of
a theorem of E.G. Effros. \textsl{Mich. Math. J.} \textbf{34}(1) (1987)
39--55.\newline
\noindent \lbrack ArhT] {A. Arhangelskii and M. Tkachenko, \textsl{%
Topological groups and related structures, }World Scientific, 2008.}\newline
\noindent \lbrack Ban] S. Banach, \textsl{Th\'{e}orie des op\'{e}rations lin%
\'{e}aires}, in: Monografie Mat., vol.1, 1932 (in: \textsl{Oeuvres}, vol.2,
PWN, 1979), translated as \textsl{Theory of linear operations}, North
Holland, 1978.\newline
\noindent \lbrack BanJ] T. Banakh and E. Jab\l o\'{n}ska, Null-finite sets
in topological groups and their applications. \textsl{Israeli J. Math}, to
appear [arXiv:1706.08155].\newline
\noindent \lbrack BanGJSJ] T. Banakh, S. G\l \k{a}b, E. Jab\l o\'{n}ska, and
J. Swaczyna, Haar-I sets: looking at small sets in Polish groups through
compact glasses, arXiv:1803.06712.

\noindent \lbrack BarJ] T. Bartoszy\'{n}ski and H. Judah, On Sierpi\'{n}ski
sets. \textsl{Proc. Amer. Math. Soc.} \textbf{108} (1990), 507-512.\newline
\noindent \lbrack BarLS] T. Bartoszy\'{n}ski, P. Larson, and S. Shelah,
Closed sets which consistently have few translations covering the line. 
\textsl{Fund. Math.} \textbf{237} (2017), 101--125.\newline
\noindent \lbrack BarS] T. Bartoszy\'{n}ski and S. Shelah, Strongly meager
sets of size continuum. \textsl{Arch. Math. Logic} \textbf{42} (2003), no.
8, 769--779.\newline
\noindent \lbrack BelS] J. Bell, A. Slomson, \textsl{Models and
ultraproducts: an introduction}, N. Holland, 1969.\newline
\noindent \lbrack BinGT] N. H. Bingham, C. M. Goldie and J. L. Teugels, 
\textsl{Regular variation}, 2nd ed., Cambridge University Press, 1989 (1st
ed. 1987). \newline
\noindent \lbrack BinO1]{\ N. H. Bingham and A. J. Ostaszewski, }Automatic
continuity: subadditivity, convexity, uniformity. \textsl{Aequationes Math.}%
\textbf{\ 78} (2009), 257--270.\newline
\noindent \lbrack BinO2]{\ N. H. Bingham and A. J. Ostaszewski, {Kingman,
category and combinatorics.} \textsl{Probability and mathematical genetics}
(Sir John Kingman Festschrift, ed. N. H. Bingham and C. M. Goldie), 135-168, 
\textsl{London Math. Soc. Lecture Notes in Mathematics} \textbf{378},
Cambridge University Press, 2010.}\newline
\noindent \lbrack BinO3]{\ N. H. Bingham and A. J. Ostaszewski, {Normed
groups: Dichotomy and duality.} \textsl{Dissert. Math. }\textbf{472} (2010),
138p.}\newline
\noindent \lbrack BinO4] N. H. Bingham and A. J. Ostaszewski, Beyond
Lebesgue and Baire II: bitopology and measure-category duality. \textsl{%
Colloq. Math.} \textbf{121} (2010), no. 2, 225--238. \newline
\noindent \lbrack BinO5] N. H. Bingham and A. J. Ostaszewski, Dichotomy and
infinite combinatorics: the theorems of Steinhaus and Ostrowski. \textsl{%
Math. Proc. Cambridge Philos. Soc.} \textbf{150} (2011), 1--22.\newline
\noindent \lbrack BinO6] N. H. Bingham and A. J. Ostaszewski,
Category-measure duality: convexity, mid-point convexity and Berz
sublinearity. \textsl{Aequationes Math.}, \textbf{91.5} (2017), 801--836 (
fuller version: arXiv1607.05750).\newline
\noindent \lbrack BinO7]{\ N. H. Bingham and A. J. Ostaszewski, {Beyond
Lebesgue and Baire IV: Density topologies and a converse Steinhaus-Weil theor%
{em}.} \textsl{Topology and its Applications} \textbf{239} (2018), 274-292
(arXiv:1607.00031).}\newline
\noindent \lbrack BinO8] N. H. Bingham and A. J. Ostaszewski, Set theory and
the analyst. \textsl{European J. Math.}, Online First, arXiv:1801.09149v2.%
\newline
\noindent \lbrack BinO9] N. H. Bingham and A. J. Ostaszewski, The
Steinhaus-Weil property: its converse, Solecki amenability and
subcontinuity, arXiv:1607.00049v3.\newline
\noindent \lbrack BinO10] {N. H. Bingham and A. J. Ostaszewski, Beyond Haar
and Cameron-Martin: topological theory, arXiv}1805.02325.\newline
\noindent \lbrack BinO11] {N. H. Bingham and A. J. Ostaszewski, }Sequential
regular variation: extensions of Kendall's theorem{, arXiv:1901.07060.}%
\newline
\noindent \lbrack Bir] G. Birkhoff, A note on topological groups. \textsl{%
Compos. Math.} \textbf{3} (1936), 427--430.\newline
\noindent \lbrack Bog] V. I. Bogachev, \textsl{Gaussian measures}, Math.
Surveys \& Monographs \textbf{62}, Amer Math Soc., 1998.\newline
\noindent \lbrack BorD] D. Borwein and S. Z. Ditor, Translates of sequences
in sets of positive measure. \textsl{Canad. Math. Bull}. \textbf{21} (1978),
no. 4, 497--498.\newline
\noindent \lbrack Bru] A. M. Bruckner, Differentiation of integrals. \textsl{%
Amer. Math. Monthly}\textbf{\ 78} (1971), no. 9, Part II, ii+51 pp.\newline
\noindent \lbrack Car] T. Carlson, Strong measure zero and strongly meager
sets. \textsl{Proc. Amer. Math. Soc.} \textbf{118} (1993), no. 2, 577--586.%
\newline
\noindent \lbrack Chr1] {J. P. R. Christensen, {On sets of Haar measure zero
in abelian Polish groups.} Proceedings of the International Symposium on
Partial Differential Equations and the Geometry of Normed Linear Spaces
(Jerusalem, 1972).\textsl{\ Israel J. Math.} \textbf{13} (1972), 255--260
(1973).}\newline
\noindent \lbrack Chr2] {J. P. R. Christensen, \textsl{Topology and Borel
structure. Descriptive topology and set theory with applications to
functional analysis and measure theory.} North-Holland Mathematics Studies 
\textbf{10}, 1974.}\newline
\noindent \lbrack CieR] K. Ciesielski, and J. Rosenblatt, Restricted
continuity and a theorem of Luzin. \textsl{Colloq. Math.} \textbf{135}
(2014), 211--225.\newline
\noindent \lbrack CrnGH] M. Crnjac, B. Gulja\v{s}, H. I. Miller, On some
questions of Ger, Grubb and Kraljevi\'{c}. \textsl{Acta Math. Hungar.} 
\textbf{57} (1991), 253--257.\newline
\noindent \lbrack Dar] {U. B. Darji, {On Haar meager sets.} \textsl{Topology
Appl.}\textbf{160} (2013), 2396--2400.}\newline
\noindent \lbrack Dev] K. Devlin, \textsl{Aspects of Constructibility},
Lecture Notes in Math. Vol. 354, Springer, 1973.\newline
\noindent \lbrack DieS] J. Diestel, A. Spalsbury,\textsl{\ The joys of Haar
measure,} Graduate Studies in Mathematics \textbf{150}, Amer. Math. Soc.,
2014.\newline
\noindent \lbrack Dra] F. Drake, \textsl{Set theory: An introduction to
large cardinals}, North-Holland, 1974.\newline
\noindent \lbrack Eff] E. G. Effros, Transformation groups and C$^{\ast }$%
-algebras. \textsl{Ann. of Math.} (2) \textbf{81} (1965), 38--55.\newline
\noindent \lbrack EleS] M. Elekes and J. Stepr\={a}ns, Less than $2^{\omega }
$ many translates of a compact nullset may cover the real line. \textsl{%
Fund. Math.} \textbf{181} (2004), 89--96.\newline
\noindent \lbrack EleT] M. Elekes and A. T\'{o}th, Covering locally compact
groups by less than $2^{\omega }$ many translates of a compact nullset. 
\textsl{Fund. Math.} \textbf{193} (2007), 243--257.\newline
\noindent \lbrack EngMS] F. van Engelen, A. W. Miller, and J. Steel, Rigid
Borel sets and better quasi-order theory. \textsl{Contemp. Math.} \textbf{65}
(1985), 199--222. \newline
\noindent \lbrack Enge] R. Engelking, \textsl{General topology}, Heldermann
Verlag, 1989.\newline
\noindent \lbrack ErdK] P. Erd\H{o}s and S. Kakutani, On a perfect set. 
\textsl{Coll. Math.} \textbf{4} (1957), 195-196.\newline
\noindent \lbrack ErdKM] P. Erd\H{o}s, K. Kunen, and D. Mauldin, Some
additive properties of sets of real numbers. \textsl{Fund. Math.} \textbf{113%
} (1981), 187--199.\newline
\noindent \lbrack FenN] J. E. Fenstad, D. Normann, On absolutely measurable
sets. \textsl{Fund. Math.} \textbf{81.2} (1973/74), 91--98.\newline
\noindent \lbrack FenMW] Qi Feng, M. Magidor, and H. Woodin,Universally
Baire sets of reals, in: \textsl{Set theory of the continuum} (Berkeley, CA,
1989), 203--242, Math. Sci. Res. Inst. Publ. \textbf{26}, Springer,1992.%
\newline
\noindent \lbrack GalMS] F. Galvin, J. Mycielski and R. M. Solovay, Strong
measure zero and infinite games. \textsl{Arch. Math. Logic} \textbf{56}
(2017), 725--732.\newline
\noindent \lbrack Hal] P. R. Halmos, \textsl{Measure Theory}, Grad. Texts in
Math. \textbf{18}, Springer 1974 (1st. ed. Van Nostrand, 1950).\newline
\noindent \lbrack HauP] O. Haupt, C. Pauc, La topologie approximative de
Denjoy envisag\'{e}e comme vraie topologie. \textsl{C. R. Acad. Sci. Paris} 
\textbf{234} (1952), 390--392.\newline
\textsl{\noindent }[Hod] W.\textsl{\ }Hodges, \textsl{A shorter model theory}%
, Cambridge University Press,1997.\newline
\noindent \lbrack Jab] E. Jab\l o\'{n}ska, Remarks on analogies between Haar
meager sets and Haar null sets, in: \textsl{Developments in functional
equations and related topics}, 149--159, Springer Optim. Appl. \textbf{124},
Springer 2017.\newline
\noindent \lbrack JasW] J. Jasi\'{n}ski and T. Weiss, Sierpi\'{n}ski sets
and strong first category. \textsl{Proc. Amer. Math. Soc.} \textbf{111}
(1991), 235-238.\newline
\noindent \lbrack Jec] T. J. Jech, \textsl{Set Theory}. 3$^{\text{rd}}$
Millennium ed. Springer, 2003.\newline
\noindent \lbrack Kak1] S. Kakutani, \"{U}ber die Metrisation der
topologischen Gruppen. \textsl{Proc. Imp. Acad. Tokyo} \textbf{12} (1936),
82--84 (reprinted in [Kak2]).\newline
\noindent \lbrack Kak2] S. Kakutani, \textsl{Selected papers.} Vol. 1. (Ed.
R. R. Kallman), Contemporary Mathematicians, Birkh\"{a}user, 1986.\newline
\noindent \lbrack Kec] A. S. Kechris: \textsl{Classical Descriptive Set
Theory.} Grad. Texts in Math. \textbf{156}, Springer, 1995.\newline
\noindent \lbrack Kel] T. Keleti, Construction of one-dimensional subsets of
the reals not containing similar copies of given patterns. \textsl{Anal. PDE}
\textbf{1} (2008), no. 1, 29--33.\newline
\noindent \lbrack Kes1] H. Kestelman, The convergent sequences belonging to
a set. \textsl{J. London Math. Soc.} \textbf{22} (1947), 130--136.\newline
\noindent \lbrack Kes2] H. Kestelman, On the functional equation $%
f(x+y)=f(x)+f(y)$. \textsl{Fund. Math.} \textbf{34}, (1947). 144--147.%
\newline
\noindent \lbrack Kuc] M. Kuczma, \textsl{An introduction to the theory of
functional equations and inequalities. Cauchy's equation and Jensen's
inequality,} 2nd ed., Birkh\"{a}user, 2009 [1st ed. PWN, Warszawa, 1985].%
\newline
\noindent \lbrack Kun] K. Kunen, \textsl{Set theory}. Studies in Logic
(London) \textbf{34}, College Publications, London, 2011.\newline
\noindent \lbrack Kur] K. Kuratowski, Ensembles projectifs et ensembles
singuliers. \textsl{Fund. Math.} 35 (1948), 131-140.\newline
\noindent \lbrack Lav1] R. Laver, On strong measure zero sets. Infinite and
finite sets (Colloq., Keszthely, 1973; dedicated to P. Erd\H{o}s on his 60th
birthday), vol. II, pp. 1025--1027. \textsl{Colloq. Math. Soc. Janos Bolyai} 
\textbf{10}, North-Holland, 1975.\newline
\noindent \lbrack Lav2] R. Laver, On the consistency of Borel's conjecture.%
\textsl{\ Acta Math.} \textbf{137} (1976), 151--169. \newline
\noindent \lbrack ManW] R. Mansfield and G. Weitkamp, \textsl{Recursive
aspects of Descriptive set theory}, Oxford Logic guides: 11, Springer, 1985.%
\newline
\noindent \lbrack Mar] N. F. G. Martin, A topology for certain measure
spaces. \textsl{Trans. Amer. Math. Soc.} \textbf{112} (1964). 1--18.\newline
\noindent \lbrack MilA1] A. W. Miller, Special subsets of the real line, in 
\textsl{Handbook of set-theoretic topology}, 201-233, N. Holland 1984.%
\newline
\noindent \lbrack MilA2] A. W. Miller, Infinite combinatorics and
definability. \textsl{Ann. Pure Appl. Logic }\textbf{41} (1989), 179--203.%
\newline
\noindent \lbrack MilM] H. I. Miller, L. Miller-Van Wieren, Some further
results on the Bowein-Ditor theorem. \textsl{Adv. Math.} \textbf{4} (2),
(2015), 121-125.\newline
\noindent \lbrack MilO] H. I. Miller and A.J. Ostaszewski, Group action and
shift-compactness. \textsl{J. Math. Anal. App.} \textbf{392} (2012), 23--39.%
\newline
\noindent \lbrack Mue] B. J. Mueller, Three results for locally compact
groups connected with Haar measure density theorem. \textsl{Proc. Amer.
Math. Soc.} \textbf{16} (1965), 1414--1416.\newline
\noindent \lbrack Ost1] A. J. Ostaszewski, Almost completeness and the
Effros open mapping principle in normed groups.\textsl{\ Topology Proc.} 
\textbf{41} (2013), 99--110.\newline
\noindent \lbrack Ost2] A. J. Ostaszewski, Beyond Lebesgue and Baire III:
Steinhhaus' Theorem and its descendants. \textsl{Topology Appl.} \textbf{160 
}(2013), 1144-1154.\newline
\noindent \lbrack Ost3] A. J. Ostaszewski, Effros, Baire, Steinhaus and
non-separability.\textsl{\ Topology Appl.} \textbf{195} (2015), 265--274.%
\newline
\noindent \lbrack Oxt] J. C. Oxtoby, \textsl{Measure and Category},
Springer, 2$^{\text{nd}}$ ed. 1980.\newline
\noindent \lbrack Pet] {B. J. Pettis, {On continuity and openness of
homomorphisms in topological groups.} \textsl{Ann. of Math. }(2) \textbf{52}
(1950), 293--308.}\newline
\noindent \lbrack Pic] {S. Piccard, {Sur les ensembles de distances des
ensembles de points d'un espace Euclidien.\ }\textsl{M\'{e}m. Univ. Neuch%
\^{a}tel} \textbf{13}, 212 pp. 1939.}\newline
\noindent \lbrack Rog1] C. A. Rogers, \textsl{Hausdorff measures}. 1$^{\text{%
st}}$ed.,1970, reprinted with a foreword by K. J. Falconer, Cambridge
University Press, Cambridge, 1998.\newline
\noindent \lbrack Rog2] {\ C. A. Rogers, J. Jayne, C. Dellacherie, F. Tops\o %
e, J. Hoffmann-J\o rgensen, D. A. Martin, A. S. Kechris, A. H. Stone, 
\textsl{Analytic sets,} Academic Press, 1980.\newline
}\noindent \lbrack Sac] Gerald E. Sacks, \textsl{Higher Recursion Theory},
Perspectives in Logic Vol 2, Springer, 1990. \newline
\noindent \lbrack Sol] {\ S. Solecki, {Amenability, free subgroups, and Haar
null sets in non-locally compact groups}. \textsl{Proc. London Math. Soc.}
(3) \textbf{93} (2006), 693--722.}\newline
\noindent \lbrack Solo] R. M. Solovay, A model of set-theory in which every
set of reals is Lebesgue measurable. \textsl{Ann. of Math.} (2) \textbf{92}
(1970), 1--56.\newline
\noindent \lbrack Ste] H. Steinhaus, Sur les distances des points de mesure
positive.\ \textsl{Fund. Math.} \textbf{1} (1920), 83-104.\newline
\noindent \lbrack Sve] R. E. Svetic, The Erd\H{o}s similarity problem: a
survey. \textsl{Real Anal. Exchange} 26 (2000/01), 525--539.\newline
\noindent \lbrack Tar1] A. Tarski, Die Wahrheitsbegriff in den
formalisierten Sprachen. \textsl{Studia Phil.} (Warsaw) \textbf{1}, 261-405
(translation in [Tar2]).\newline
\noindent \lbrack Tar2] A. Tarski, The concept of truth in formalized
languages. In: \textsl{Logic, semantics, metamathematics, papers from
1923-1938} (translated by J. H. Woodger), 152-278, Oxford, 1956.

\noindent \lbrack TopH] {F. Tops\o e and J. Hoffmann-J\o rgensen, {Analytic
spaces and their application}, in [Rog2, Part 3].}\newline
\noindent \lbrack Vid] Z. Vidny\'{a}nszky, Transfinite inductions producing
coanalytic sets. \textsl{Fund. Math.} \textbf{224} (2014), 155--174.\newline
\noindent \lbrack Wei] A. Weil,\ \textsl{L'int\'{e}gration dans les groupes
topologiques}, Actualit\'{e}s Scientifiques et Industrielles 1145, Hermann,
1965 (1$^{\text{st }}$\ ed. 1940).\newline

\bigskip

\noindent Faculty of Engineering and Natural Sciences/Mathematics,
International University of Sarajevo, 71000 Sarajevo, Bosnia-Herzegovina;
harrymiller609@yahoo.com\newline
Faculty of Engineering and Natural Sciences/Mathematics, International
University of Sarajevo, 71000 Sarajevo, Bosnia-Herzegovina;
lmiller@ius.edu.ba\newline
Mathematics Department, London School of Economics, Houghton Street, London
WC2A 2AE; A.J.Ostaszewski@lse.ac.uk\newpage

\end{document}